\documentclass[twocolumn]{autart}    

\usepackage{graphicx}          

 \usepackage{amsfonts,amsmath,amssymb,bm}
\usepackage{color,cases,float}
\usepackage{graphicx}
\usepackage{hyperref}
\usepackage{psfrag}
\usepackage{eucal}
\usepackage{algorithm}
\usepackage{overpic}
\usepackage{subfigure}
\usepackage{url}
\usepackage[all,cmtip]{xy}
\usepackage{algorithm}
\usepackage[noend]{algpseudocode}

\definecolor{darkblue}{rgb}{0.0,0.0,0.6}
\hypersetup{colorlinks,breaklinks,
	linkcolor=darkblue,urlcolor=darkblue,
	anchorcolor=darkblue,citecolor=darkblue}
\newtheorem{assumption}{Assumption}
\newtheorem{definition}{Definition}
\newtheorem{corollary}{Corollary}
\newtheorem{theorem}{Theorem}
\newtheorem{remark}{Remark}

\makeatletter
\newcommand{\setassumptiontag}[1]{
  \let\oldtheassumption\theassumption
  \renewcommand{\theassumption}{#1}
  \g@addto@macro\endassumption{
    \addtocounter{assumption}{-1}
    \global\let\theassumption\oldtheassumption}
  }
\makeatother

\renewcommand{\qed}{\hfill\blacksquare}

\newcommand{\E}{\mathrm{E}}

\newcommand{\R}{\mathrm{Re}}
\newcommand{\Proj}{\mathrm{Proj}}

\def\blab{\boldsymbol{\lambda}}
\def\fb{\mathbf{f}}
\def\R{\mathbb{R}}
\def\Wb{\mathbf{W}}

\def\zb{\mathbf{z}}
\def\bx{\boldsymbol{x}}

\def\ab{\boldsymbol{a}}

\def\mb{\boldsymbol{m}}

\def\s{\sigma}
\def\g{\gamma}

\def\bmu{\boldsymbol{\mu}}
\def\hmu{\hat{\mu}}
\def\hbmu{\hat{\boldsymbol{\mu}}}

\def\kb{\boldsymbol{k}}
\def\lb{\boldsymbol{l}}

\def\gb{\boldsymbol{g}}
\def\fb{\boldsymbol{f}}
\def\R{\mathbb{R}}
\def\Wb{\boldsymbol{W}}

\def\zb{\mathbf{z}}
\def\xb{\boldsymbol{x}}

\def\bmu{\boldsymbol{\mu}}

\def\Mb{\boldsymbol{M}}
\def\Qb{\boldsymbol{Q}}
\def\bla{\boldsymbol{\lambda}}

\def\Rb{\boldsymbol{R}}
\def\Sb{\boldsymbol{S}}
\def\Ab{\boldsymbol{A}}
\def\Pb{\boldsymbol{P}}

\def\bxi{\boldsymbol{\xi}}
\def\r{~}
\def\BibTeX{{\rm B\kern-.05em{\sc i\kern-.025em b}\kern-.08em
		T\kern-.1667em\lower.7ex\hbox{E}\kern-.125emX}}

\allowdisplaybreaks

\begin{document}

\begin{frontmatter}

\title{Convergence Rate of Generalized Nash Equilibrium Learning in Strongly Monotone Games with Linear Constraints} 

\thanks[footnoteinfo]{This paper was not presented at any IFAC 
meeting. Corresponding author T. Tatarenko. Tel. +49 6151 16-25048. 
}

\author[Paestum]{Tatiana Tatarenko}\ead{tatiana.tatarenko@tu-darmstadt.de},    
\author[Rome]{Maryam Kamgarpour}\ead{maryam.kamgarpour@epfl.ch},               

\address[Paestum]{TU Darmstadt, Germany}  
\address[Rome]{EPFL, Switzerland}             

\begin{keyword}                           
Multi-agent systems, Stochastic control and game theory, Convex optimization                          
\end{keyword}                             

\begin{abstract}                          
We consider payoff-based learning of a generalized Nash equilibrium (GNE)  in multi-agent systems. Our focus is on  games with jointly convex constraints of a linear structure  and strongly monotone pseudo-gradients. We present a convergent procedure based on a partial regularization technique and establish the convergence rate of its iterates under one- and two-point payoff-based feedback. 
To the best of our knowledge, this work is the first one characterizing the convergence speed of iterates to a variational GNE in the class of games under consideration.  
\end{abstract}

\end{frontmatter}

\section{Introduction}\label{sec:intro}

Generalized Nash equilibrium (GNE) problems model settings where multiple decision makers have coupled objective and constraint functions. The applicability of GNE problems spans from economics to engineering systems  \cite{facchinei2007generalized}.  Several works have focused on characterizing and computing a solution to the GNE problem given knowledge of the cost and constraint functions, see \cite{jordan2023first} for a recent overview of algorithms.  Our focus here is on computing a solution to a class of GNE problems without prior knowledge on players' cost or constraint functions. We assume only access to the so-called ``payoff" information. Namely,  players can evaluate their cost and constraint functions at a  chosen joint action. Hence, our goal is to compute a GNE solution based on these function evaluations.


Learning equilibria under payoff-based information setting is a well-established field. It is motivated  by the fact that in many practical multi-agent interactive settings, the agents do not know functional form of the objectives and constraints and can only access the values of their objective functions at a played action. Such situations arise for example, in electricity markets (unknown price functions or constraints), network routing (unknown traffic demands/constraints), and sensor coverage problems (unknown density function on the mission space). Even in the cases where cost and constraint functions are known, a payoff-based approach can be  advantageous due to the complexity of computing gradients. This latter motivation  is in line with works on derivative-free  zeroth-order optimization of complex functions.

Most past work on computing a solution to a GNE problem assumes availability of gradients. The general approach to solve this problem class in convex games is to cast it as an equivalent quasi-variational inequality (QVI) problem. In the case in which all agents share the same coupling constraint, a subset of GNE solutions, referred to as variational GNE (v-GNE),  can be found through a variational inequality (VI) problem. Moreover, a v-GNE is shown to be a refinement of the GNE and to possess economic merit~\cite{KULKARNI2012}.  

The focus of early algorithms for GNE problems has been on proving \textit{asymptotic} convergence rather than characterizing \textit{convergence rates}~\cite{facchinei2007generalized}. Recently, \cite{jordan2023first,meng2023linear} proposed algorithms for GNE problems with monotone pseudo-gradients and  linear coupling constraints, and established the convergence rates of their respective algorithms. The work of \cite{meng2023linear}  considered the case where all players share the same constraint function and derived a primal-dual approach with a linear convergence rate. This work assumed  a strongly monotone pseudo-gradient, and the full row rank of the constraint's matrix, based on which a specific notion of convergence was defined. On the other hand, \cite{jordan2023first} considered the more general class of GNE problems where agents can have individually different coupling constraints, derived an  algorithm converging to an approximate Karush-Kuhn-Tucker (KKT) solution  and established its convergence rate.
However, the notion of approximate solution in the work above requires checking an approximate KKT condition at each iteration and generally does not imply convergence of the iterates to a GNE solution. Furthermore,  both of the above works required knowledge of the constraints and the costs' gradients.

Considering the payoff-based information setting, for games without coupling constraints  past work had established convergence rates  under strong monotonicity of the pseudo-gradient \cite{drusvyatskiy2022improved} or strong variational stability of equilibria \cite{TatKamECC24}. For games with jointly convex coupling constraints, 
\cite{tat_kam_TAC}  derived a payoff-based convergent algorithm under the assumption that the game admits a potential function. The potential function assumption  simplifies the analysis as it enables one to cast the problem as   constrained optimization. Thus, \cite{suli} extended this setting to include non-potential games with jointly linear constraints, and presented a Tikhonov regularization approach, similar to \cite{tatarenko2019learning}. However, \cite{suli} did not establish a convergence rate for the algorithm.

To the best of our knowledge, deriving convergence rate of payoff-based learning of a variational Nash equilibrium in a GNE problem has not been addressed. Thus, our paper is devoted to addressing this problem.  

We develop a convergent payoff-based algorithm for GNE problems with linear coupling constraints under the standard assumptions of strongly monotone and Lipschitz continuous pseudo-gradients. 
To obtain the convergence rate in our setting  we introduce a regularization to the VI formulation of our problem, and establish bounds between solutions of the regularized problem and that of the original problem, see Lemmas~\ref{lem:reg_t},~\ref{lem:epsApprox}, and ~\ref{lem:t_vs_t-1}.  

Preliminary results related to our approach were developed in~\cite{ECC25}. In contrast to the present article, \cite{ECC25} considered 2-point feedback and ignored local agents' constraints. The 1-point feedback leads to an unbounded variance of the gradient estimators, which needs to be taken into account by an appropriate choice of the algorithm's parameters. Furthermore, local constraints challenges the randomized gradient estimation techniques. 
In summary, our contributions are as follows:
\begin{itemize}
\item We present a convergent payoff-based  algorithm for jointly convex GNE problems with strongly monotone pseudo-gradients and linear coupling constraints;
\item We establish the convergence rate of the algorithm's iterates as $O(1/t^{1/4-\delta})$ under the 1-point  and $O(1/t^{1/2-\delta})$ under the 2-point feedback, where $\delta$ is an arbitrary small parameter;  
\item We demonstrate that these rates can be improved to $O(1/t^{4/15-\delta})$ and $O(1/t^{4/7})$, respectively, in the case of an interior v-GNE.
\end{itemize}

\subsection*{Notations. }The set $\{1,\ldots,N\}$ is denoted by $[N]$. We consider real normed space $\R^D$ and denote by $\R_{+}^D$ its non-negative halfspace. The column vector $\xb\in\R^D$ is denoted by $\xb = [x^1,\ldots, x^D]$. We use superscripts to denoted coordinates of vectors and the player-related functions and sets. 
For any function $f:K\to\R$, $K\subseteq\R^D$, $\nabla_{x^i} f(\xb) = \frac{\partial f(\xb)}{\partial x^i}$ is the partial derivative taken in respect to the $x^i$th variable (coordinate) in the vector argument $\xb\in\R^D$.
 We use $\langle \cdot,\cdot\rangle$ to denote the inner product in $\R^D$.
We use $\|\cdot\|$ to denote the Euclidean norm induced by the standard dot product in $\R^D$. 
A mapping $g:\R^D\to \R^D$ is said to be \emph{strongly monotone} on $Q\subseteq \R^D$ with the constant $\nu$, if for any $u, v\in Q$, $\langle g(u)-g(v), u - v\rangle \ge \nu\|u - v\|^2$; \emph{strictly monotone}, if the strict inequality holds for $\nu=0$ and $u\neq v$, and \emph{merely monotone} if the inequality holds for $\nu=0$.
We use $\Proj_{\Omega}{v}$ to denote the projection of $v\in \R^D$ onto a set $\Omega\subseteq \R^D$.
The  expectation of a random value $\xi$ is denoted by $\E\{\xi\}$. Its conditional expectation with respect to some $\sigma$-algebra $\EuScript F$ is denoted by $\E\{\xi|\EuScript F\}$.
We use the big-$O$ notation, that is, the function $f(x): \R\to\R$ is $O(g(x))$ as $x\to a$ for some $a\in\R$, i.e. $f(x)$ = $O(g(x))$ as $x\to a$, if $\lim_{x\to a}\frac{|f(x)|}{|g(x)|}\le K$ for some positive constant $K$.

\section{Game Settings}
We consider a non-cooperative game $\Gamma$ between $N$ players. Let player $i$, $i \in [N]$, choose its local action from $\R^{d_i}$ and $J^i: \R^D\to\R$ with $D=\sum_{i=1}^Nd_i$ denote her cost function. We assume existence of a local constraint set $\Ab^i\subseteq \R^{d_i}$ for each player $i\in[N]$ as well as of a \emph{shared coupling constraint set} $C$ defined by the following system of linear inequalities: 
\begin{align}\label{eq:coupled_cs}
C = \{\ab\in\R^D: \gb(\ab)=K\ab-\lb\le \boldsymbol{0}\},
\end{align}
where $K\in\R^{n\times D}$ is a matrix with columns   $\kb^{i,j} \in \R^n$, for $j\in [d_i]$, $i\in[N]$, $\lb\in\R^{n}$, $\ab = [\ab^1,\ldots,\ab^N]\in\R^D$, with $\ab^i = [a^{i,1},\ldots,a^{i,d_i}]\in\R^{d_i}$ being the action of  player $i$ and $\ab^{-i}$ denoting the joint action of all players except player $i$. 
Thus, the feasible action set in game $\Gamma$ is $C\cap\Ab$, where $\Ab = \Ab^1\times\ldots\times\Ab^N$, and the coupled action set of the agent $i$ is $\Ab^i\cap{C}^i(\ab^{-i})$, where ${C}^i(\ab^{-i}) = \{\ab^i\in\R^{d_i}: K\ab-\lb\le \boldsymbol{0}\}$. 

A \emph{generalized Nash equilibrium} (GNE) in  game $\Gamma$ represents a joint action from which no player has any incentive to unilaterally deviate.
\begin{definition}\label{def:GNE}
 A joint action $\ab^*\in  C\cap\Ab$ in the game $\Gamma=\Gamma(N, \{\Ab^i\}, \{J^i\}, C)$ with the coupling constraints $C$ is called a \emph{generalized Nash equilibrium} (GNE) if $J^i(\ab^{*,i},\ab^{*,-i})\le J^i(\ab^i,\ab^{*,-i})$ for any $i\in[N]$ and $\ab^i\in\Ab^i\cap C^i(\ab^{*,-i})$.
 
Observe that if $C=\R^{D}$ then $ C^i(\ab^{-i})=\{\ab^i: \ab^i\in \Ab^{i}\}$ and any $\ab^*$ for which the above inequality holds is a \emph{Nash equilibrium (NE)} in the game with no coupling constraints.
 \end{definition}

We consider convex games as defined below. 
\begin{assumption}\label{assum:convex}
 The game under consideration is \emph{convex}. Namely, for all $i\in[N]$ the cost function $J^i(\ab^i, \ab^{-i})$ is defined on $\R^{D}$, continuously differentiable in $\ab$ and convex in $\ab^i$ for  fixed $\ab^{-i}$. 
\end{assumption}
Given the differentiability of the cost functions, we can define the game's pseudo-gradient mapping below.
\begin{definition}\label{def:mapping}
The mapping $\Mb = [\Mb^1, \ldots, \Mb^N]:\R^{D}\to\R^{D}$, referred to as  the \emph{pseudo-gradient} of  $\Gamma(N, \{\Ab^i\}, \{J^i\}, C)$ is defined by
 \begin{align}\label{eq:gamemapping}
 &\Mb^i(\ab) = [M^{i,1}(\ab), \ldots, M^{1,d_i}(\ab)]^{\top}, \cr
 &M^{i,j}(\ab)= \frac{\partial J^i(\ab)}{\partial a^{i,j}}, \; \ab\in  \Ab\cap C, j\in[d_i], i\in[N].
 \end{align}
\end{definition}
It is known that characterizing and computing GNE for game $\Gamma$ is challenging \cite{facchinei2007generalized}. However, a solution concept referred to as  \emph{variational generalized Nash equilibrium} can be characterized through a variational inequality defined below and will be the focus of our work  \cite{facchinei2007generalized}. 
\begin{definition}\label{def:VGNE}
    The set of \emph{variational generalized Nash equilibria  (v-GNE)} in the game $\Gamma$ is defined as the set of solutions to the following variational inequality problem:
    \begin{align}\label{eq:GVI}
    \nonumber
    \mbox{Find } &\ab^*\in  \Ab\cap C: \, \langle\Mb(\ab^*),\ab-\ab^*\rangle\ge 0, \\
    &\mbox{ for all }\ab\in\Ab\cap C .\tag{$VI_{(\Mb, \Ab\cap C)}$}
\end{align}
\end{definition}
The next lemma states the relation between GNE and v-GNE (see \cite{facchinei2007generalized} for the proof). 
\begin{lem}\label{lem:GNEandVGNE}
    Under Assumption\r\ref{assum:convex}, the set of generalized Nash equilibria in  game $\Gamma$ contains the set of variational generalized Nash equilibria, i.e. v-GNE $\subseteq$ GNE\footnote{Without loss of generality we assume that for any set $\EuScript S$ the following holds: $\emptyset\subseteq \EuScript S$, where $\emptyset$ denotes an empty set.}.   
\end{lem}



\subsection{Formulation as a game without coupling constraint}
To deal with the coupling constraints in the game $\Gamma$, we follow the approach from \cite{tat_kam_TAC} and define an uncoupled game $\Gamma^{\blab}$, for which we introduce a virtual player, the so-called \emph{dual} player indexed by $N+1$. The action of the player $N+1$ is denoted by  $\blab\in\Ab^{N+1}=\R^{n}_{+}$, and corresponds to the dual variable of the coupling constraint. 
Thus, given  shared constraints defined by a linear set $C = \{\ab: \gb(\ab)\le \boldsymbol{0}\}$, the augmented game $\Gamma^{\blab}$ contains $N+1$ players, the cost function $U^{i}$ of each of the first $N$ players, referred to as primal players,  is defined as her local Lagrangian: 
\begin{align}\label{eq:J^{pr}r}
U^{i}(\ab,\blab) = J^i(\ab) + \langle \blab, \gb(\ab)\rangle.
\end{align}
The cost function of the dual player $N+1$ is defined as
\begin{align}\label{eq:J_du}
U^{N+1}(\ab,\blab) = -  \langle \blab, \gb(\ab)\rangle.
\end{align}
The joint action set in the uncoupled game $\Gamma^{\blab}$ is denoted by $\Ab^{\blab} = \Ab^1\times\ldots\times \Ab^N\times\Ab^{N+1}$. Hence, $\Gamma^{\blab} = \Gamma^{\blab}(N+1, \{U^{i}\}, \Ab^{\blab})$ is a game with uncoupled constraints. We denote a joint action of this game by $[\ab,\blab]$ and refer to $\ab\in\Ab$ and $\blab\in\Ab^{N+1}=\R^{n}_{+}$ as  its \emph{primal and dual part} respectively.  

In this work,  we focus on  games with strongly monotone pseudo-gradients to be able to establish uniqueness of v-GNE and convergence rate of our algorithm to it. 
\begin{assumption}\label{assum:strmon}
The pseudo-gradient $\Mb$ is \emph{strongly monotone} over $\R^D$, i.e. there exists a constant $\nu>0$ such that $\langle\Mb(\ab_1)-\Mb(\ab_2),\ab_1-\ab_2\rangle\ge\nu\|\ab_1-\ab_2\|^2$ for any $\ab_1,\ab_2\in\R^D$.
\end{assumption}
We are ready to formulate the result connecting solutions of the initial game $\Gamma$ with that of game $\Gamma^{\blab}$.
\begin{lem}\label{lem:exist_uncoupled}
 Let $\Gamma(N,\{\Ab^i\}, \{J^i\}, C)$ be a game with shared constraints defined by a convex set $C = \{\ab: \, \gb(\ab)\le 0\}$. 
 Then,
 
 1) given Assumption~\ref{assum:convex},  if $[\ab^*,\bla^*]$ is an NE in $\Gamma^{\blab}$, then $\ab^*$ is a v-GNE of $\Gamma$. Moreover, if $\ab^*$ is a v-GNE of $\Gamma$, then there exists $\blab^*\in\R^n_{+}$ such that $[\ab^*,\bla^*]$ is an NE in $\Gamma^{\blab}$,
 
 2) given Assumptions~\ref{assum:convex} and~\ref{assum:strmon}, there exists a NE in $\Gamma^{\blab}$,
 
 3) given Assumption~\ref{assum:strmon}, if $[\ab_1^*,\bla_1^*]$ and $[\ab_2^*,\blab_2^*]$ are two NE in $\Gamma^{\blab}$, then $\ab_1^* = \ab_2^*$,
 
 4) for any NE $[\ab^*, \bla^*]$ in $\Gamma^{\blab}$ there exists a constant $\Lambda > 0$ such that $\|\bla^*\|\le \Lambda$.
 \end{lem}
Please refer to Appendix~\ref{app:exist_uncoupled} for the proof.
The lemma above implies that one can aim to find a solution (NE) to the game $\Gamma^{\blab}$ with uncoupled action sets while seeking a v-GNE in the initial game $\Gamma$. By solving $\Gamma^{\blab}$, one achieves a v-GNE. Moreover, under Assumption~\ref{assum:strmon}, as stated above, this v-GNE exists and is unique. 

\subsection{Uncoupled game properties}
Our goal here is to derive properties of the uncoupled game, which we will utilize in establishing convergence rate of our payoff-based algorithm.

Consider the  specific settings of the affine shared constraints,  $\gb(\ab) = K\ab-\lb$.
In this setting the mapping $\Wb = [\Wb^{1},\ldots,\Wb^{N},\Wb^{N+1}]:\R^{D+n}\to\R^{D+n}$, denoting  the \emph{pseudo-gradient} of the uncoupled game $\Gamma^{\blab}$, is given by
 \begin{align}\label{eq:gamemapping_aug}
  &\Wb^{i}(\ab, \blab)=[W^{i,1}(\ab, \blab), \ldots, W^{i,d_i}(\ab, \blab)], \cr
 &\qquad\mbox{where }W^{i,j}(\ab,\blab)= M^{i,j}(\ab) + \frac{\partial \langle\bla, K\ab-\lb\rangle}{\partial a^{i,j}}\cr
 &= M^{i,j}(\ab) + \langle \blab,\kb^{i,j}\rangle ,\; j\in[d_i], \quad i\in[N],\cr
 &\Wb^{N+1}(\ab, \blab)= -K\ab+\lb.
 \end{align}
Let us denote the primal part of the mapping $\Wb$ by 
\begin{align}\label{eq:primalM}
\Wb^{pr}(\ab,\blab) = [\Wb^{1}(\ab,\blab),\ldots,\Wb^{N}(\ab,\blab)]\in\R^D.
\end{align}

     

Notice that according to ~\ref{eq:gamemapping_aug}, the pseudo-gradient $\Wb$ of the uncoupled game $\Gamma^{\blab}$ is linear in the dual action $\blab$. Thus, $\Gamma^{\blab}$ is a convex game and its Nash equilibria  are characterized by solutions to the following variational inequality:  
\begin{align}\label{eq:VI1}
    \mbox{Find }&[\ab^*,\blab^*]\in\Ab^{\blab}: \, \langle\Wb(\ab^*,\blab^*),[\ab,\blab]-[\ab^*,\blab^*]\rangle\ge 0, \cr&\,\mbox{ for all }[\ab,\blab]\in\Ab^{\blab} \tag{$VI$}.
\end{align} 
Our goal is to solve the variational inequality~\eqref{eq:VI1} so as to obtain a v-GNE of the original game. The challenge lies in the fact that $\Wb$ is not strongly monotone. To address this, we suitably regularize the mapping and connect the solutions of the regularized problem to that of the original problem $\Gamma^{\bla}$.

Consider the  Tikhonov regularization  with time-varying term $\varepsilon_t I_{0,n}$ introduced to the mapping $\Wb$, where $I_{0,n}\in\R^{D+n\times D+n}$ is a diagonal matrix with the first $D$ zero-elements and the last $n$ to be equal to $1$. Thus, the following regularized pseudo-gradient is constructed: 
\begin{align}\label{eq:reg}
     \Wb_t(\ab,\blab) = \Wb(\ab,\blab) + [\underbrace{0,\ldots, 0}_{D}, \varepsilon_t\blab].
\end{align}

\begin{lem}
\label{lem:reg_t}
Under Assumptions~\ref{assum:convex}, \ref{assum:strmon}, the regularized mapping $\Wb_t(\ab,\blab)$ is strongly monotone, namely,  for any $ [\ab_1,\blab_1]\in\Ab^{\blab}$, $[\ab_2,\blab_2]\in\Ab^{\blab}$
\begin{align}\label{eq:strMon_W_t}
\langle\Wb_t(\ab_1,\blab_1)&-\Wb_t(\ab_2,\blab_2),[\ab_1,\blab_1]-[\ab_2,\blab_2]\rangle \cr&\ge \nu\|\ab_1-\ab_2\|^2 + \varepsilon_t\|\blab_1-\blab_2\|^2.
\end{align}
\end{lem}

\begin{pf}
Due to the \emph{linear structure of the coupling constraints}, it is straightforward to verify that the following relation holds for any two actions $ [\ab_1,\blab_1]\in\Ab^{\blab}$, $[\ab_2,\blab_2]\in\Ab^{\blab}$ of the uncoupled game $\Gamma^{\blab}$.
\begin{align}\label{eq:aug-psedo}
    \langle\Wb(\ab_1,\blab_1)-\Wb(\ab_2,\blab_2),[\ab_1,\blab_1]&-[\ab_2,\blab_2]\rangle \cr&\ge \nu\|\ab_1-\ab_2\|^2.
\end{align}
Thus, the statement of the lemma holds due to Assumption~\ref{assum:strmon} and definition of the mapping $\Wb_t$ in~\eqref{eq:reg}. 
    
$\qed$\end{pf}
Given the result above one can use a gradient-descent type approach to solve the following regularized variational inequality:   Find $[\ab_t^*,\blab_t^*]\in\Ab^{\blab}$ such that
\begin{align}\label{eq:VI2} 
 \langle\Wb_t(\ab_t^*,\blab_t^*),[\ab,\blab]-[\ab_t^*,\blab_t^*]\rangle\ge0, \,\forall [\ab,\blab]\in\Ab^{\blab} \tag{$VI_t$}.
\end{align}
The solution of the above is unique due to Lemma \ref{lem:reg_t}.

The lemma below estimates the  rate with which the primal part of the  solution to the regularized~\eqref{eq:VI2} converges to the primal part of the solution to the variational inequality~\eqref{eq:VI1} as $\varepsilon_t\to 0$.
\begin{lem}\label{lem:epsApprox}
    Let the mapping $\Wb(\ab,\blab)$ be as in Definition~\ref{eq:gamemapping_aug} with $\Mb(\ab)$ satisfying Assumption~\ref{assum:strmon}. 
    Let $[\ab^*,\blab^*]$ be a minimal norm solution to~\eqref{eq:VI1} and $[\ab^*_t,\blab^*_t]$ be a solution to~\eqref{eq:VI2}. Then, as $\varepsilon_t\to 0$, $[\ab^*_t,\blab^*_t]$ converges to $[\ab^*,\blab^*]$ and, moreover,
    \[\|\ab^*-\ab_t^*\|^2 \le\varepsilon_t\frac{\|\blab^*\|\|\blab^*-\blab^*_t\|}{\nu}\le\varepsilon_t\frac{\|\blab^*\|^2}{\nu}.\]
\end{lem}
See Appendix~\ref{app1} for the proof.

The next lemma bounds the distance between the solutions to two consecutive regularized ~\eqref{eq:VI2}'s. 
\begin{lem}\label{lem:t_vs_t-1}
Consider  $\Wb(\ab,\blab)$ in~\eqref{eq:gamemapping_aug} with $\Mb(\ab)$ satisfying Assumption~\ref{assum:strmon}. 
    Let $[\ab^*_t,\blab^*_t]$ and $[\ab^*_{t-1},\blab^*_{t-1}]$ be solutions to~\eqref{eq:VI2} and $(VI_{t-1})$ respectively. Then the following relation holds: 
    \begin{align*}
        \|\ab^*_t-\ab^*_{t-1}\|^2 &= O\left(\frac{(\varepsilon_t-\varepsilon_{t-1})^2}{\varepsilon_{t}}\right),\cr
        \|\blab^*_t-\blab^*_{t-1}\|^2 &= O\left(\frac{(\varepsilon_t-\varepsilon_{t-1})^2}{\varepsilon^2_{t}}\right)
    \end{align*}
    
\end{lem}
\begin{pf}
    See Appendix~\ref{app1} for the proof.
$\qed$\end{pf}

The above three lemmas on the properties of \eqref{eq:VI2} and the relation of its solutions to that of \eqref{eq:VI1} will be used in the proof of algorithm convergence in Section \ref{sec:convergence}. 

\section{Payoff-based Learning algorithm}\label{sec:procedure}

The \emph{payoff-based} settings assume the following information. Each primal player $i$ in the initial game $\Gamma$ has access to the value of her cost function $J^i(t) = J^i(\ab(t))$ at every time step $t$, given the corresponding joint action $\ab(t)$. Moreover, the value of the constraint function, $\gb(t)=K\ab(t)-\lb$, is revealed to each player $i\in[N]$ at every $t$. In this setting, neither the full information on the joint action $\ab(t)$ nor the closed form of the cost and shared constraint functions are available to the agents.


We start setting up the learning algorithm based on the payoff-based information described above.

\setlength{\textfloatsep}{15pt}
\begin{algorithm}[t!]
	\caption{One-point and two-point payoff-based algorithm for learning  Nash equilibria }\label{alg:algorithm1}
	\begin{algorithmic}[!t]
		\Require Action set $\Ab^i \subset \R^{d_i}$, the sequences $\{\sigma_t\},  \{\gamma_t\}$, $\{\varepsilon_t\}$ initial states $\bmu^i(0)$, $\blab(0)$, $i\in[N]$.
		\For {$t = 0,1, \ldots$}
		\State  Sample a query point $\bxi^i(t)$ according to probability density 
               \eqref{eq:density}. Play action $\ab^i(t) = \Proj_{\Ab^i}[\bxi^i(t)]$ and observe $\gb(t) = K\ab(t)-\lb$, $J^i(t) = J^i(\ab(t))$.
		  
	     	\noindent{\color{blue} { \footnotesize\ttfamily  /* Simultaneously and similarly,  other players sample their  query points $\bxi^{-i}(t)$.   /* }}
		\State One-point scheme: Perform the one-point gradient estimate $\mb_c^i(t)$ according to \eqref{eq:est_Gd}.		
           \State Two-point scheme:  Observe additional  $U^{i}_0(t)$ (see~\eqref{eq:J_0}) and perform the two-point gradient estimate according to \eqref{eq:est_Gd2}.
		\State Update the state according to \eqref{eq:alg}:
  	  \State$\qquad\bmu^i(t+1)=\Proj_{\Ab^i}[\bmu^i(t)-\gamma_t\mb^i_c(t)], \, 
                   i\in[N]\mbox{ (primal players)};$
           \State$\qquad\blab(t+1)=\Proj_{\R^{n}_{+}}[\blab(t)-\gamma_t(-\gb(t)+\varepsilon_t\blab(t))], \,\mbox{ (dual player)};$
		\noindent{\color{blue} { \footnotesize\ttfamily  /* Simultaneously and similarly,  other players update their states (primal and dual).  /* }}
		\EndFor
		
	\end{algorithmic}
\end{algorithm}

\subsection{Algorithm iterates}
For each time $t$ in the game $\Gamma^{\bla}$, let us denote the actions of the primal player $i$ and the dual player $N+1$  by $\ab^i(t)$ and $\blab(t)$ respectively. 

We define $\mb^i_c(t)$, $c\in\{1,2\}$ to be an estimate of the pseudo-gradient's elements of the primal players in the game $\Gamma^{\bla}$ at time $t$. In words, $\mb^i_c(t)$, $c\in\{1,2\}$ estimates  the vector  $\Wb^i(\ab(t),\blab(t))$, $i\in[N]$, see Definition ~\eqref{eq:gamemapping_aug}. Here, $c=1$ denotes the one-point and $c=2$ denotes the two-point procedure estimate case, respectively. The exact estimation procedure will be detailed in Subsection \ref{sub:grad_estimate}.
 
The algorithm iterates corresponding to the primal player $i$, $i\in[N]$, and the dual player $N+1$ are denoted by  $\bmu^i$  and  $\blab$, respectively. These iterates are updated  as follows:
\begin{align}
	\label{eq:alg}
	&\bmu^i(t+1)=\Proj_{\Ab^i}[\bmu^i(t)-\gamma_t\mb^i_c(t)],\cr
        &\blab(t+1)=\Proj_{\R^{n}_{+}}[\blab(t)-\gamma_t(-\gb(t)+\varepsilon_t\blab(t))],
   \end{align}
where $\bmu^i(0)\in \R^{d_i}$, $\blab(0)\in\R^{n}_{+}$ are  arbitrary finite vectors, $\g_t$ is the step size or the learning rate,  $\gb(t)$ is the value of the constraint function available to the player $i\in[N]$ at time $t$, i.e. $\gb(t)=K\ab(t)-\lb$, and the choice of $\ab(t)$ will be described below.

Note that according to the discussion in the previous section, the procedure in~\eqref{eq:alg} mimics the gradient descent play in  game $\Gamma^{\blab}$. The regularization term $\varepsilon_t\blab(t)$ is added only in the update of the dual players. The step size $\gamma_t$ should be chosen based on the regularization parameter of the game, as well as on the bias and variance of the pseudo-gradient estimate $\mb_c^i(t)$. The term $\mb_c^i(t)$ is obtained using payoff-based feedback as described below.

\subsection{Gradient estimations}
\label{sub:grad_estimate}
To estimate the unknown gradients of each player's objective with respect to her own decision variable, and using only payoff-based information, we use the well-established technique of randomized sampling \cite{spall1997one}. In our setting, we use the Gaussian distribution for sampling inspired by \cite{Thatha,NesterovSpokoiny}.

Given $\bmu^i(t)\in\Ab^i$, consider the iterates, projected on the shrunk set $(1-\rho_t)\Ab^i$ for some $\rho_t\ge 0$, 
\begin{align}\label{eq:muhat}
    \hat{\bmu}^i(t) = \Proj_{(1-\rho_t)\Ab^i}\bmu^i(t), \quad i\in[N]
\end{align}
Each player $i$, $i\in[N]$, samples the random vector $\bxi^i(t)$ according to the multidimensional normal distribution $\EuScript N(\hbmu(t) = [\hmu^{i,1}(t),\ldots,\hmu^{i,d_i}(t)]^{\top},\sigma_t)$ with the following density function:
\begin{align}\label{eq:density}
	p^i&(\bx^i;\hbmu^i(t),\sigma_{t})\cr=& \frac{1}{(\sqrt{2\pi}\sigma_{t})^{d_i}}\exp\left\{-\sum_{k=1}^{d_i}\frac{(x^{i,k}-\hmu^{i,k}(t))^2}{2\sigma^2_{t}}\right\},
\end{align}
To query the objectives at a locally feasible action, each player projects her action on the local constraint sets, namely, $\Proj_{\Ab^i}\bxi^i(t)$. Thus, the joint action is
\begin{align}\label{eq:a}
\ab(t) = \Proj_{\Ab}\bxi(t)=[\Proj_{\Ab^1}\bxi^1(t),\ldots,\Proj_{\Ab^N}\bxi^N(t)]^{\top}.
\end{align}
The evaluated cost function at the above point, in the game $\Gamma^{\blab}$ at $\ab(t)$ and $\bla(t)$ (see definition~\eqref{eq:J^{pr}r}) is denoted by 
\[U^{i}(t) =U^{i}(\ab(t),\blab(t))=J^i(\ab(t)) + \langle \blab(t), \gb(\ab(t)),\]
is revealed to each primal player $i\in[N]$ in $\Gamma^{\blab}$. 

Note that the value $U^{i}(t)$ is available to player $i$, $i=1\in[N]$ as follows. First, according to the payoff-based information settings,  agent $i$ has access to the values $J^i(t) = J^i(\ab(t))$ and $\gb(t)=\gb(\ab(t))$. Second, given $\blab^i(0) = \blab(0)$ for any player $i$, namely, same initial dual vector $\blab(0)$ for any player $i$, the update in Equation \eqref{eq:alg} can be run locally by each primal player and the common dual variable $\blab(t)$ will be available for each such player at each time $t$. Alternatively, the dual vector $\blab(t)$ can be  broadcast by a central coordinator, playing the role of a virtual dual player. Note that an appropriate coordination of  the dual variable value $\blab(t)$ is required and commonly considered \cite{tat_kam_TAC,suli} to guarantee convergence of the algorithm to a v-GNE.

After describing the randomization for choosing action queries above, we can discuss the gradient estimators. 

In the one-point setting ($c=1$), player $i$  estimates her local pseudo-gradient $\Wb^i(\cdot,\cdot)$ evaluated at the point of the joint iterate $(\hbmu(t),\blab(t))$ with $\hbmu(t)=[\hbmu^1(t),\ldots,\hbmu^N(t)]^{\top}$  as follows:
\begin{align}\label{eq:est_Gd}
	\mb^i_1(t) = { U^{i}(t)}\frac{{\bxi^i(t)} -\hbmu^i(t)}{\sigma^2_t}.
\end{align}
In the two-point setting ($c=2$) at each iteration $t$, each player $i$ makes two queries: a query corresponding to the chosen action $\ab^i(t)$, $i\in[N]$, and another query of the cost function $J^i$ and the constraint function $g^i$ at the point
$\hat{\bmu}(t) = \Proj_{(1-\rho_t)\Ab}\bmu(t)$, $i\in[N]$.
Hence, there is an extra piece of information available to each player, namely the cost function value at the mean vector:
\begin{align}\label{eq:J_0}
U^{i}_0(t) = U^{i}(\hbmu(t),\blab(t)) = J^i(\hbmu(t)) + \langle \blab(t), \gb(\hbmu(t)). 
\end{align}
Then each player uses the following estimation of the local pseudo-gradient $\Mb^i(\cdot)$ at the point  $\hbmu(t)$:
	\begin{align}\label{eq:est_Gd2}
		\mb^i_2(t) = (U^{i}(t) - U^{i}_0(t))\frac{{\bxi^i(t)} -\hbmu^i(t)}{\sigma^2_t}.
	\end{align}
 Given the estimation approaches for the cases $c=1,2$ above, the procedure's iterates are summarized in Algorithm~\ref{alg:algorithm1}.

 \begin{remark}
    For the randomzied sampling technique in gradient estimation, one can also use the uniform distribution over the unit sphere.  The resulting gradient estimates will have the  same properties as in the case of the Gaussian distribution choice (see, for example, \cite{bravo2018bandit}). 
\end{remark}

\begin{remark}
    As it can be expected, our analysis will demonstrate a better behavior of the algorithm for $c=2$ corresponding to the two-point feedback.
    This is due to the bounded variance of the estimate $\mb_2^i(t)$ for any choice of the parameter $\sigma_t$.   
    Alternatively, to guarantee a bounded variance of the gradients' estimates, one can use memory, namely the value of the  function $U^i(t-1)$ at the previous iterate $t-1$, instead of $U^i_0(t)$, in order to obtain $\mb_2^i(t)$ as it has been proposed in~\cite{Yuanhanqing24}. 
 \end{remark}

\begin{remark}
    Algorithm~\ref{alg:algorithm1} uses three projection operators in relation to the local action sets $\{\Ab^i\}_i$. They are in the states' updates and actions' generation as described by Equations ~\eqref{eq:alg},~\eqref{eq:muhat} and~\eqref{eq:a} respectively. Thus, we assume the local action sets $\Ab^i$ to be ``projection-friendly" (the Euclidean projection onto them can be done with a small computational effort). Moreover, observe that in Equation~\eqref{eq:muhat},  the projection of a feasible point $\bmu(t)\in\Ab$ is onto the shrunk set  $(1-\rho_t)\Ab$. The projection onto the shrunk set  $(1-\rho_t)\Ab$ can be replaced by shifting the point $\bmu(t)$ along a feasible direction. This idea is presented, for example, in~\cite{bravo2018bandit} and is called feasibility adjustment. 
    

    
\end{remark}

\subsection{Interpretation as stochastic gradient descent}

We provide insight into the procedure defined in~\eqref{eq:alg}  by deriving an analogy to a  stochastic gradient algorithm and we use this insight in proving convergence rate of the algorithm. Denote 
\begin{align}\label{eq:densityfull}
p( \bx; \hbmu, \sigma)=\prod_{i=1}^{N}p^i(\bx^i;\hbmu^i,\sigma)
\end{align}
as the  density function of the joint distribution of the players' joint query point $\bxi$, given some mean vector $\hbmu =[\hbmu^1,\ldots,\hbmu^N]$ (see~\eqref{eq:density} for the definition of the individual density function $p^i$). For any $\sigma > 0$, $i\in[N]$, and given any $\blab\in\R^{n}_{+}$ define $ \tilde{U}^{i}_{\sigma}(\hbmu,\blab) : \R^{D} \rightarrow \R$ as
\begin{align}
\label{eq:mixedJ}
\tilde{U}^{i}_{\sigma} &(\hbmu,\blab)= \int_{\mathbb R^{D}}U^{i}(\bx,\blab)p( \bx; \hbmu, \sigma)d\bx.
\end{align}
Thus, $\tilde{U}^{i}_{\sigma}$, $i\in[N]$, is the $i$th primal player's cost function in the mixed strategies of the game $\Gamma^{\blab}$, where the strategies are sampled from the Gaussian distribution with the density function in~\eqref{eq:densityfull}.
For $i\in[N]$ define $\tilde{\Wb}^{i}_{\sigma} (\cdot)=[\tilde W^{i,1}_{\sigma}(\hbmu), \ldots, \tilde W^{i,d_i}_{\sigma}(\hbmu)]^{\top}$
as the $d_i$-dimensional mapping with the following elements:
\begin{align}\label{eq:mapp2}
\tilde W^{i,k}_{\sigma} (\hbmu,\blab)=\frac{\partial {\tilde U^{i}_{\sigma}(\hbmu,\blab)}}{\partial \mu^{i,k}}, \mbox{ for $k\in[d_i]$}.
\end{align}
Furthermore, let us define 
\begin{align}\label{eq:hatmu}
   \hat{\mb}_c^{i}(t) =\begin{cases}
	U^{i}(\bxi(t),\blab(t))\frac{\bxi^i(t) -\hbmu^i(t)}{\sigma^2_t}, \, \mbox{if $c=1$},\\
	(U^{i}(\bxi(t),\blab(t)) - U^{i}_0(t))\frac{\bxi^i(t) -\hbmu^i(t)}{\sigma^2_t}, \mbox{ if $c=2$}.
	\end{cases} 
\end{align}
Thus, the term
\[{\mb}_c^{i}(t) = \hat{\mb}_c^{i}(t) +\Pb^{i}(t),\]
where
\begin{align}\label{eq:Pterm}
&\Pb^{i}(t) =  (U^{i}(t) - U^{i}(\bxi(t),\blab(t))) \frac{\bxi^i(t) -\hbmu^i(t)}{\sigma^2_t} \cr
&=  (U^{i}(\ab(t),\blab(t)) - U^{i}(\bxi(t),\blab(t))) \frac{\bxi^i(t) -\hbmu^i(t)}{\sigma^2_t}
\end{align}
can be interpreted as the term caused by taking projection of the sampled point $\bxi(t)$ onto the feasible set $\Ab$ to obtain the action $\ab(t)$. 
Next, let $\Rb_{c}^{i}(t)$, $c=1,2$, denote:
\begin{align}\label{eq:R}
&\Rb_{c}^{i}(t) = \hat{\mb}_c^{i}(t) - \tilde{\Wb}^{i}_{\sigma_t} (\hbmu(t),\blab(t)). 
\end{align}
Finally, we introduce the term \[\Qb^i(t)=\tilde\Wb^{i}_{\sigma_t}(\hbmu(t),\blab(t))-\Wb^{i}(\hbmu(t),\blab(t)),\] which evaluates the bias between the primal player $i$th pseudo-gradient in the game $\Gamma^{\blab}$ and its counterpart in the mixed strategies, $i\in[N]$. With the above definitions, the joint update rule presented in \eqref{eq:alg}  is equivalent to:
\begin{align}
\label{eq:pbavmuQ}
\bmu(t+1) =\Proj_{\Ab}[\bmu(t)& -\gamma_t\big({\Wb}^{pr}(\hbmu(t),\blab(t))\cr
&+\Qb(t)+\Pb(t)+\Rb_c(t)\big)].
\end{align}
Recall that the cases $c=1$, $c=2$ above correspond to the one-point and two-point gradient estimators, respectively. Moreover, above we use the notation $\Wb^{pr}(\cdot,\cdot)$ defined in~\eqref{eq:primalM} as well as  the notations $\bmu(t) = [\bmu^1(t),\ldots,\bmu^N(t)]$, $\Qb(t) = [\Qb^1(t),\ldots,\Qb^N(t)]$, $\Pb(t) = [\Pb^1(t),\ldots,\Pb^N(t)]$, and $\Rb_c(t) = [\Rb_c^1(t),\ldots,\Rb_c^N(t)]$.
In the next sections devoted to the convergence analysis of Algorithm~\ref{alg:algorithm1} we will use the representation~\eqref{eq:pbavmuQ} for the updates of the primal-related variables comprised by $\bmu(t)$, whereas, according to~\eqref{eq:alg}, the dual-related action, namely ${\blab}(t)$, is updated as follows: 
\begin{align}\label{eq:dual}
	\blab(t+1)=\Proj_{\R^{n}_{+}}[\blab(t)-\gamma_t(-K\hat\bmu(t) &+ \Sb(t) + \lb\cr
 &+\varepsilon_t\blab(t))],
   \end{align}
where $\Sb(t) =  K(\hat\bmu(t) - \ab(t))$.


To be able to conduct the analysis of the algorithm, we make the following assumptions which are standard in the literature on  generalized Nash equilibria seeking \cite{jordan2023first}.
\begin{assumption}\label{assum:compact}
 For all $i\in[N]$ the set $\Ab^i$ is \emph{compact}. 
\end{assumption}

\begin{assumption}	\label{assum:infty}
	Each function $J^{i}(\bx) = O(\|\bx\|^{2})$ as $\|\bx\|\to\infty$.
\end{assumption}
\begin{remark}\label{rem:Lip}
    Assumption~\ref{assum:infty} above is required, as we will consider a payoff-based approach to v-GNE learning, and to this end, will estimate gradients based on randomized sampling, namely on a 
    a distribution with an unbounded support. Moreover, this assumption will play an important role in proving boundedness of the iterates in the learning procedure.
    Observe that  Assumption~\ref{assum:infty} is equivalent to the commonly made assumption that the functions $\Mb^i$, $i\in[N]$, are \emph{Lipschitz continuous on $\R^{D}$} with some constant $L>0$. 
\end{remark}

\subsection{Properties of the  gradient estimators}
We now proceed with formulation of the main properties of the terms $\Rb_c(t)$, $c=1,2$, $\Qb(t)$, $\Pb(t)$, and $\Sb(t)$. 
Let $\EuScript F_{t}$ be the $\sigma$-algebra generated by the random variables $\{\bmu(k),\blab(k),\bxi(k)\}_{k\le t}$. The next lemmas characterize the stochastic term $\Rb^i_c(t)$ and the bias term $\Qb^i(t)$, $i\in[N]$. In particular, they state that the conditional expectation of $\Rb^i_c(t)$  is equal to zero, whereas its variance depends on the settings of $\sigma_t$ in the case $c=1$ and is uniformly bounded in the case $c=2$. As for the term $\Qb^i(t)$, its variance is proportional to the choice of the parameter $\sigma_t$. The proof of the following lemmas (Lemmas~\ref{lem:sample_grad}-\ref{lem:Qsq}) can be found in Lemma~1 and~3 in \cite{tatarenko2019learning}.
\begin{lem}\label{lem:sample_grad}
Given Assumptions\r\ref{assum:convex} and~\ref{assum:infty}, for $c=1,2$, almost surely 
$\E\{\Rb_{c}^{i}(t)|\EuScript F_t\} = 0$ for all $i\in[N]$.
\end{lem}
\begin{lem}\label{lem:Rsq}
Under Assumptions~\ref{assum:convex} and~\ref{assum:infty},  for $c=1,2$,  given $\lim_{t\to\infty}\sigma_t = 0$, 
\[\E\{\|\Rb_{c}^{i}(t)\|^2 | \EuScript F_t\} = \begin{cases}
	O\left(\frac{\|\blab(t)\|^2}{\sigma_t^2}\right), \, &\mbox{ if $c=1$},\\
	O(\|\blab(t)\|^2), \, &\mbox{ if $c=2$}.
\end{cases}\]
\end{lem}
Finally, the following relation holds for the second moment of the term $\Qb(t)$. 
\begin{lem}\label{lem:Qsq}
Under Assumption~\ref{assum:infty}, as $\s_t\to 0$, \\
$\E\{\|\Qb^{i}(t)\|^2| \EuScript F_t\} = O(\s^2_t).$
\end{lem}
Next, we characterize the second moment of the terms $\Pb(t)$ in~\eqref{eq:pbavmuQ} and $\Sb(t)$ in~\eqref{eq:dual}.
\begin{lem}\label{lem:Psq}
Let Assumptions~\ref{assum:convex} and~\ref{assum:infty} hold. Choose $\sigma_t, \rho_t$ such that $\lim_{t\to\infty}\sigma_t = 0, \lim_{t\to\infty}\rho_t = 0, \lim_{t\to\infty}\frac{\rho_t}{\s_t}=\infty $. Then almost surely
\begin{align*}
   \E\{\|\Pb^{i}(t)\|^2 | \EuScript F_t\}  =O\left((1+\|\blab(t)\|^2)\left(\frac{e^{-\frac{{\rho}_t^2}{2\sigma_t^2}}}{\sigma_t^{D+2}}\right)\right).
\end{align*}
\end{lem}
\begin{pf} 
 According to the definition of the term $\Pb^i(t)$ in~\eqref{eq:Pterm} and the function $U^{i}(\ab(t),\blab(t))$ in~\eqref{eq:J^{pr}r},
 \begin{align}\label{eq:Psq}
 \|\Pb^i(t)\|^2 &\le 2(|J^i(\ab(t))-J^i(\bxi(t))|^2\cr
 &\,+ \|\blab(t)\|^2\|\gb(\ab(t))-\gb(\bxi(t))\|^2)\frac{\|\bxi^i(t) -\hat\bmu^i(t)\|^2}{\sigma^4_t}\cr
 &=2(|J^i(\Proj_{\Ab}\bxi(t))-J^i(\bxi(t))|^2 \cr
 &\,+ \|\blab(t)\|^2\|\gb(\Proj_{\Ab}\bxi(t))-\gb(\bxi(t))\|^2)\cr
 &\qquad\qquad\qquad\qquad\times\frac{\|\bxi^i(t) -\hat\bmu^i(t)\|^2}{\sigma^4_t}.
 \end{align}
 Next, we can follow the line of proof in Proposition 1 3) in~\cite{TatKamTCNS24}  to conclude that for any vector-function $\fb = (f^1,\ldots, f^k)$ with $f^j:\R^D\to \R$ satisfying Assumptions~\ref{assum:convex} and~\ref{assum:infty} the following holds:
 \begin{align}\label{eq:projterm1}
\E\{\| \fb(\Proj_{\Ab}\bxi(t))-\fb(\bxi(t)) \|^2|\EuScript F_t\} = O\left(\frac{e^{-\frac{{\rho}_t^2}{2\sigma_t^2}}}{\sigma_t^{D}}\right).
 \end{align}
 Combining~\eqref{eq:projterm1} with~\eqref{eq:Psq} and taking into account linearity of $\gb$, we obtain the result.
$\qed$\end{pf}


\begin{lem}\label{lem:Sterm2moment} Under Assumptions~\ref{assum:convex} and~\ref{assum:infty}, almost surely
\begin{align*}
    \E\{\|\Sb(t)\|^2 | \EuScript F_t\}  =O\left(\sigma_t^2+\frac{e^{-\frac{{\rho}_t^2}{2\sigma_t^2}}}{\sigma_t^{D+2}}\right).
\end{align*}
\end{lem}
\begin{pf}
We notice that 
\begin{align*}
\Sb(t) &= K(\hat\bmu(t) - \ab(t)) \cr
&= K(\hat\bmu(t) - \bxi(t) +\bxi(t) - \Proj_{\Ab}\bxi(t)),
\end{align*} 
which implies
\begin{align*}
\|\Sb(t)\|^2 = O(\|\hat\bmu(t) - \bxi(t)\|^2)   + O(\|\bxi(t) - \Proj_{\Ab}\bxi(t)\|^2) 
\end{align*} 
Taking conditional expectation of both sides and applying~\eqref{eq:projterm1} as well as $\E\{O(\|\hat\bmu(t) - \bxi(t)\|^2| \EuScript F_t\} = O(\sigma_t^2)$, we conclude the result.
$\qed$\end{pf}
The above properties enable us to establish the main results as follows. 

\section{Convergence rate of the algorithm}
\label{sec:convergence}
In this section we provide the main result on convergence of Algorithm~\ref{alg:algorithm1} and its convergence  rate, see Theorem \ref{th:main} and Corollary \ref{cor:opt_rate}, respectively. Furthermore, in Subsection \ref{sub:th_intNE} we show that if the unique v-GNE in the game $\Gamma$ belong to the interior of the set $\Ab$, the convergence rate can be improved, see Theorem \ref{th:th_intNE} and Corollary \ref{cor:int}. 


\begin{theorem}\label{th:main}
Let Assumptions~\ref{assum:convex}-\ref{assum:infty} hold.  Let the joint action $\ab^*$ be the   unique v-GNE of the game $\Gamma$. 
Let  the parameters in Algorithm~\ref{alg:algorithm1} be chosen  as follows:
\begin{align*}
    &\gamma_t = \frac{G_c}{t^{g_c}}, \, \varepsilon_t=\frac{E_c}{t^{e_c}}, \, 
      \sigma_t=\frac{S_c}{t^{s_c}}, \, \rho_t =\frac{R_c}{t^{r_c}},
      \, \mbox{ for $c=1, 2$},
\end{align*}
where $g_c,e_c,s_c,r_c>0$ with $s_c>r_c$, $g_c+e_c<1$.
Let $h_1 = \min\{2-g_1-e_1,  g_1+r_1, 2g_1-2s_1\}$, $h_2 = \min\{2-g_2-e_2,  g_2+r_2, 2g_2\}$ be such that $h_c-g_c>0$.

Then for the iterates in Algorithm~\ref{alg:algorithm1} the following holds: 
\begin{align*}
    \E\|\bmu(t)-\ab^*\|^2 =O\left(\frac{1}{t^{\min\{e_c,h_c-g_c\}}}\right), \, \mbox{for $c=1, 2$}.
\end{align*}
\end{theorem}
\begin{pf}
Existence and uniqueness of a v-GNE in the game $\Gamma$ follows from Lemma~\ref{lem:exist_uncoupled} 2)-3).
To prove the above rate, we apply a four-step convergence analysis as follows. First, we prove boundedness of the iterates $\blab(t)$.  In the second step,  using the above boundedness, we  estimate the distance between the  vector $\zb(t)=[\bmu(t),\blab(t)]$ and a regularized solution $\zb_t^*= [\ab^*_t, \blab^*_t]$. The third step consists in using the rate with which $\E\|\blab(t) - \blab_t^*\|$ tends to zero to specify the rate of convergence for the primal part, namely for the term $\E\|\ab(t) - \ab_t^*\|$. Finally, the fourth step uses the result provided in Lemma~\ref{lem:epsApprox} to conclude the convergence rate in expectation.


To proceed with the above steps, we leverage the algorithm's representation~\eqref{eq:pbavmuQ}-\eqref{eq:dual} and use the results provided in the previous subsection on its terms $\Qb(t)$, $\Pb(t)$, $\Rb_c(t)$, and $\Sb(t)$.  We also introduce $\hat\zb(t)$: $\hat\zb(t) = [\hat \bmu(t),\blab(t)]$.

\textbf{Step 1. Boundedness of iterates.} Taking into account properties  of the terms $\Pb(t)$, $\Qb(t)$, $\Sb(t)$, and $\Rb_c(t)$ (see Lemmas~\ref{lem:sample_grad}-\ref{lem:Qsq}),  we conclude that the procedure~\eqref{eq:pbavmuQ}-\eqref{eq:dual} can be considered as a perturbed gradient-based method to solve the variational inequality $(VI_t)$. Thus, in the following we will estimate the distance between the algorithm's iterate at time $t+1$ and the solution $\zb^*_t=[\ab^*_t,\blab^*_t]$ to $(VI_t)$. We obtain
\begin{align}\label{eq:ineq1}
&\|\zb(t+1)-\zb_t^*\|^2\le \|\bmu(t) -\gamma_t\big(\Wb^{pr}(\hat\zb(t))+\Qb(t)\cr
&\qquad+\Pb(t)+\Rb_c(t)) - \ab_t^*\|^2  \cr
&+\|\blab(t)-\gamma_t(-K\hat\bmu(t) + \Sb(t) + \lb+\varepsilon_t\blab(t)) - \blab_t^*\|^2\cr
\le&2\gamma_t\langle\Wb_t(\zb_t^*),\zb(t) - \zb_t^*\rangle+\|\zb(t) - \zb_t^*\|^2\cr
&-2\gamma_t\langle\Wb^{pr}(\hat\zb(t)), \bmu(t)-\ab^*_t\rangle \cr
&-2\gamma_t\langle-K\hat\bmu(t) +  \lb + \varepsilon_t\blab(t),\blab(t)-\blab^*_t\rangle\cr 
& -2\gamma_t\langle\Qb(t)+\Pb(t)+\Rb_c(t),\bmu(t) - \ab_t^*\rangle\cr
&+\gamma_t^2\|\Wb^{pr}(\hat\zb(t))+\Qb(t)+\Pb(t)+\Rb_c(t)\|^2\cr
&-2\gamma_t\langle\Sb(t),\blab(t)-\blab^*_t\rangle \cr 
&\quad+\gamma_t^2\|-K\hat\bmu(t) + \Sb(t) + \lb+\varepsilon_t\blab(t)\|^2,
\raisetag{1.35\normalbaselineskip}
\end{align}
where the first inequality is due to the non-expansive property of the projection operator, whereas in the second inequality we used the fact that $\zb_t^*$ solves $(VI_t)$ and, thus, $\langle\Wb_t(\zb_t^*),\zb(t) - \zb_t^*\rangle\ge 0$. 
We notice that according to the definition of the mappings $\Wb^{pr}$ and $\Wb_t$  (see~\eqref{eq:primalM} and \eqref{eq:reg}) the following holds for the second and the third terms of the right hand side in the last inequality:
\begin{align}\label{eq:eq1_1}
    -2\g_t(\langle&\Wb^{pr}(\hat\zb(t)),\bmu(t)-\ab^*_t\rangle\cr
    &\quad+ \langle-K\hat\bmu(t) +  \lb + \varepsilon_t\blab(t),\blab(t)-\blab^*_t\rangle) \cr
    &=-2\g_t\langle\Wb_t(\hat\zb(t)),\zb(t) - \zb_t^*\rangle\cr
    &= -2\g_t(\langle\Wb_t(\zb(t)),\zb(t) - \zb_t^*\rangle \cr
    &\quad+ \langle\Wb_t(\hat\zb(t))-\Wb_t(\zb(t)),\zb(t) - \zb_t^*\rangle).
\end{align}
Above, in the last equality we added and subtracted $-2\g_t\langle\Wb_t(\zb(t)),\zb(t) - \zb_t^*\rangle$. 
Next, we estimate the term $\langle\Wb_t(\hat\zb(t))-\Wb_t(\zb(t)),\zb(t) - \zb_t^*\rangle$ in the equality above.
According to the definition of $\Wb_t$ (see~\eqref{eq:reg}), we rewrite this term as 
\begin{align}\label{eq:ineq1_1}
    &\langle\Wb_t(\hat\zb(t))-\Wb_t(\zb(t)),\zb(t) - \zb_t^*\rangle \cr
    &= \langle\Wb^{pr}(\hat\zb(t)) - \Wb^{pr}(\zb(t)),\bmu(t)-\ab^*_t\rangle \cr
    &\qquad+ \langle K(\bmu(t)-\hat\bmu(t)), \blab(t)-\blab_t^*\rangle\cr
    &=\langle\Mb(\hat\bmu(t))-\Mb(\bmu(t)),\bmu(t)-\ab^*_t\rangle\cr
    &\qquad+\langle K(\bmu(t)-\hat\bmu(t)), \blab(t)-\blab_t^*\rangle\cr
   & \le \|\Mb(\hbmu(t))-\Mb(\bmu(t))\|\|\bmu(t) - \ab_t^*\| \cr
   &\qquad+ \|K\|\|\bmu(t) - \hat\bmu(t)\|\|\blab(t)-\blab_t^*\|\cr
    &\le L\|\hbmu(t)-\bmu(t)\|\|\bmu(t) - \ab_t^*\| \cr
    &\qquad+ \|K\|\|\bmu(t) - \hat\bmu(t)\|\|\blab(t)-\blab_t^*\| \cr
    &= L\rho_t\|\bmu(t) - \ab_t^*\| + \|K\|\rho_t\|\blab(t)-\blab_t^*\|.
\end{align}
In the above, the first inequality is due to the Cauchy-Schwarz one, the second inequality is due to Assumption~\ref{assum:infty}, where $L$ is the Lipschitz constant of $\Mb$, see Remark~\ref{rem:Lip}, and the last equality is due to the fact that $\|\hbmu(t)-\bmu(t)\| = \|\Proj_{(1-\rho_t)\Ab}\bmu(t)-\bmu(t)\| =O(\rho_t)$, see Proposition 3 in~\cite{TatKamTCNS24}.
  Next, observe that Lemma~\ref{eq:strMon_W_t} implies 
 \begin{align*}
    \langle&\Wb_t(\zb(t))-\Wb_t(\zb_t^*),\zb(t) - \zb_t^*\rangle \cr&\ge \nu\|\bmu(t)-\ab_t^*\|^2 + \varepsilon_t\|\blab(t)-\blab_t^*\|.
\end{align*}
Plugging this last inequality together with~\eqref{eq:ineq1_1} and~\eqref{eq:eq1_1} into~\eqref{eq:ineq1}, we conclude that 
\begin{align}\label{eq:ineq2}
\|&\zb(t+1)-\zb_t^*\|^2\le(1-2\gamma_t\nu)\|\bmu(t)-\ab_t^*\|^2 \cr
&\qquad+ (1-2\gamma_t\varepsilon_t)\|\blab(t)-\blab_t^*\|^2\cr
&\qquad+2\gamma_tL\rho_t\|\bmu(t) - \ab_t^*\| + 2\gamma_t\|K\|\rho_t\|\blab(t)-\blab_t^*\|\cr
&\qquad-\gamma_t\langle\Qb(t)+\Pb(t)+\Rb_c(t),\bmu(t) - \ab_t^*\rangle\cr
&\qquad+\gamma_t^2\|\Wb^{pr}(\hat{\zb}(t))+\Qb(t)+\Pb(t)+\Rb_c(t)\|^2\cr
&\qquad-2\gamma_t\langle\Sb(t),\blab(t)-\blab^*_t\rangle\cr
&\qquad+\gamma_t^2\|-K\hat\bmu(t) + \Sb(t) + \lb+\varepsilon_t\blab(t)\|^2.
\end{align}
Next,  we take the conditional expectation of the both sides in the last inequality with respect to $\EuScript F_t$. Thus, we obtain
\begin{align}\label{eq:ineq3}
 \E&\{\|\zb(t+1)-\zb_t^*\|^2|\EuScript F_t\}
\le(1-2\gamma_t\nu)\|\bmu(t)-\ab_t^*\|^2 \cr
&\qquad+ (1-2\gamma_t\varepsilon_t)\|\blab(t)-\blab_t^*\|^2 + O(\gamma_t\rho_t) \cr
&\qquad+ 2\gamma_t\|K\|\rho_t\|\blab(t)-\blab^*\|\cr
&\qquad+2\gamma_t\E\{(\|\Qb(t)\|+\|\Pb(t)\|)\|\bmu(t) - \ab_t^*\| | \EuScript F_t\}\cr
&\qquad+2\gamma_t\E\{\langle\Rb_c(t),\bmu(t) - \ab_t^*\rangle | \EuScript F_t\}\cr
 &\qquad+2\gamma_t^2\E\{(2\|\Wb^{pr}(\hat{\zb}(t))\|^2+2\|\Qb(t)\|^2\cr
 &\qquad\qquad+2\|\Pb(t)\|^2+2\|\Rb_c(t)\|^2) | \EuScript F_t\}\cr
&\qquad+2\gamma_t\E\{\|\Sb(t)\|\|\blab(t)-\blab^*_t\| | \EuScript F_t\} \cr
&\qquad+ 2\gamma_t^2\E\{(2\|K\|^2\|\hbmu(t)\|^2+2\|\Sb(t)\|^2+2\|\lb\|^2\cr
&\qquad\qquad+2\varepsilon_t^2\|\blab(t)\|^2) | \EuScript F_t\}\cr
&\le(1-2\gamma_t\nu)\|\bmu(t)-\ab_t^*\|^2 + (1-2\gamma_t\varepsilon_t)\|\blab(t)-\blab_t^*\|^2\cr
&+O\left(\g_t^2+\frac{\g_t^2}{\s_t^2}\mathbb{I}_{c=1}\right)\|\blab(t)-\blab_t^*\|^2+O\left(\gamma_t\rho_t+\g_t^2\right)\cr
&\qquad\qquad\qquad\times(\|\bmu(t)-\ab_t^*\| + \|\blab(t)-\blab_t^*\|),
\raisetag{1.8\normalbaselineskip}
\end{align}
 where $\mathbb{I}_{c=1}$ is the indicator function equal to 1, if $c=1$, and equal to 0, if $c=2$. The first inequality above is due to the Cauchy-Schwarz inequality. In the second inequality we use Lemma~\ref{lem:sample_grad} implying that $\E\{\langle\Rb_c(t),\bmu(t) - \ab_t^*\rangle | \EuScript F_t\} = \boldsymbol{0}$, the definition of the mapping $\Wb^{pr}$, the fact that the vectors $\bmu(t)$, $\hat\bmu(t)$, $\ab_t^*$ have bounded norms (belong to the compact set $\Ab$), and Lemmas~\ref{lem:Rsq}-\ref{lem:Qsq} bounding the moments of the terms $\Pb(t)$, $\Qb(t)$, and $\Sb(t)$.
 Moreover, we used the setting of the parameters $\gamma_t$, $\sigma_t$, $\varepsilon_t$, and $\rho_t$ to obtain the term of the big-O order. 
 Next, we combine~\eqref{eq:ineq3} with the following relations: 
 \begin{align}\label{eq:rel}
     &\|\bmu(t)-\ab_t^*\|\le \|\bmu(t)-\ab_{t-1}^*\| +\|\ab_{t}^* -\ab_{t-1}^*\|,\cr
     &\|\bmu(t)-\ab_{t-1}^*\| \le \|\bmu(t)-\ab_{t-1}^*\|^2 + 1,\cr
     &\|\blab(t)-\blab_t^*\|\le\|\blab(t)-\blab_{t-1}^*\| + \|\blab_{t-1}-\blab_t^*\|,\cr
     &\|\blab(t)-\blab_{t-1}^*\| \le  \|\blab(t)-\blab_{t-1}^*\|^2 + 1,\cr
     &\|\bmu(t)-\ab_t^*\|^2\le (1+0.5\gamma_t\nu)\|\bmu(t)-\ab_{t-1}^*\|^2\cr
     &\qquad\qquad\qquad+ \left(1+\frac{2}{\gamma_t\nu}\right)\|\ab_t^*-\ab_{t-1}^*\|^2,\cr
     &\|\blab(t)-\blab_t^*\|^2\le (1+0.5\gamma_t\varepsilon_t)\|\blab(t)-\blab_{t-1}^*\|^2 \cr
     &\qquad\qquad\qquad+ \left(1+\frac{2}{\gamma_t\varepsilon_t}\right)\|\blab^*_t-\blab_{t-1}^*\|^2,\cr
     &\|\ab^*_t-\ab^*_{t-1}\|^2 = O\left(\frac{(\varepsilon_t-\varepsilon_{t-1})^2}{\varepsilon_{t}}\right), \cr
     &\|\blab^*_t-\blab^*_{t-1}\|^2 = O\left(\frac{(\varepsilon_t-\varepsilon_{t-1})^2}{\varepsilon^2_{t}}\right) \mbox{(see Lem.~\ref{lem:t_vs_t-1})}
 \end{align}
 to conclude that for sufficiently large $t$,
 \begin{align}\label{eq:2}
 & \E\{\|\zb(t+1)-\zb_t^*\|^2|\EuScript F_t\}\le(1-\gamma_t\nu)\E\|\bmu(t)-\ab_{t-1}^*\|^2 \cr
 &\qquad\qquad\qquad+ (1-\gamma_t\varepsilon_t)\E\|\blab(t)-\blab_{t-1}^*\|^2\cr
 &+O\left(\frac{(\varepsilon_t-\varepsilon_{t-1})^2}{\gamma_t\varepsilon^3_{t}}+\gamma_t\rho_t+\g_t^2+\frac{\g_t^2}{\s_t^2}\mathbb{I}_{c=1}\right)\cr
 &\times(\underbrace{\|\bmu(t)-\ab_{t-1}^*\|^2 + \|\blab(t)-\blab_{t-1}^*\|^2}_{\|\zb(t)-\zb_{t-1}^*\|^2}+1).
 \raisetag{3.1\normalbaselineskip}
\end{align}
We notice that the projection of the iterates for the dual variable $\blab(t)$, is taken onto unbounded set $\R^{n}_{+}$. To be able to get the convergence result, we follow the idea of Lemma~3 in \cite{tatarenko2019learning}, to demonstrate almost sure boundedness of the procedure's iterates.  For this purpose we first notice that, given the settings of the time-dependent parameters, we have
\[\sum_{t=1}^{\infty}\frac{(\varepsilon_t-\varepsilon_{t-1})^2}{\gamma_t\varepsilon^3_{t}}+\gamma_t\rho_t+\g_t^2+\frac{\g_t^2}{\s_t^2}\mathbb{I}_{c=1}<\infty.\]
With the inequality above in place, we recall the result from  Theorem 2.5.2 in \cite{NH} (see Appendix~\ref{app_bound} of this paper and let $X(t) = \zb(t)$, $\xb=\zb$, $V(t,\xb) = \|\xb-\zb_{t-1}^*\|^2$ in Theorem~\ref{th:bound}) to conclude the uniform boundedness of the norm $\|\zb(t)-\zb_{t-1}^*\|$ almost surely and, thus, uniform boundedness of $\|\blab(t)\|$ for any $t$ almost surely.
\textbf{Step 2. Convergence rate of $\E\|\blab(t+1)-\blab_{t}^*\|^2$.} 
By taking full expectation of both sides in~\eqref{eq:2}, we get for sufficiently large $t$ and both $c=1,2$:
\begin{align}\label{eq:c=1}
&\E\|\zb(t+1)-\zb_t^*\|^2\le\left(1-\frac{G_c\nu}{t^{g_c}}\right)\E\|\bmu(t)-\ab_{t-1}^*\|^2 \cr
&\qquad+ \left(1-\frac{G_cE_c}{t^{g_c+e_c}}\right)\E\|\blab(t)-\blab_{t-1}^*\|^2+O\left(\frac{1}{t^{h_c}}\right)\cr
&\le\left(1-\frac{G_cE_c}{t^{g_c+e_c}}\right)\E\|\zb(t)-\zb_{t-1}^*\|^2+O\left(\frac{1}{t^{h_c}}\right),
\raisetag{1.8\normalbaselineskip}
\end{align}
where $h_1 = \min\{2-g_1-e_1,  g_1+2r_1,  2g_1-2s_1\}$ and $h_2 = \min\{2-g_2-e_2,  g_2+r_2, 2g_2\}$. 
Thus, according to the Chung's lemma (see Lemma~\ref{lem:chung} in Appendix~\ref{app2}), we conclude that $\E\|\zb(t)-\zb^*_{t-1}\|^2 = O(1/t^{h_c-g_c-e_c})$. Thus, 
$$\E\|\blab(t+1)-\blab_{t}^*\|^2 = O(1/t^{h_c-g_c-e_c}).$$
\textbf{Step 3. Convergence rate of $\E\|\bmu(t)-\ab_{t-1}^*\|^2$.} Using this last relation together with the relation $\E\|\zb(t+1)-\zb_t^*\|^2\ge \E\|\bmu(t+1)-\ab_{t}^*\|^2+\left(1-\frac{G_cE_c}{t^{g_c+e_c}}\right)\E\|\blab(t+1)-\blab_{t}^*\|^2$ in the first inequality of~\eqref{eq:c=1}, we get
\begin{align}\label{eq:c=1_1}
    &\E\|\bmu(t+1)-\ab_{t}^*\|^2\le\left(1-\frac{G_c\nu}{t^{g_c}}\right)\E\|\bmu(t)-\ab_{t-1}^*\|^2 \cr
    &+ \left(1-\frac{G_cE_c}{t^{g_c+e_c}}\right)\cr
    &\qquad\times\left( O\left(\frac{1}{t^{h_c-g_c-e_c}}\right)-O\left(\frac{1}{(t+1)^{h_c-g_c-e_c}}\right)\right)\cr
&\qquad\qquad\qquad+O\left(\frac{1}{t^{h_c}}\right)\cr
&= \left(1-\frac{G_c\nu}{t^{g_c}}\right)\E\|\bmu(t)-\ab_{t-1}^*\|^2 +O\left(\frac{1}{t^{h_c}}\right),
\raisetag{1.8\normalbaselineskip}
\end{align}
since $O\left(\frac{1}{t^{h_c-g_c-e_c}}\right)-O\left(\frac{1}{(t+1)^{h_c-g_c-e_c}}\right) = O\left(\frac{1}{t^{1+h_c-g_c-e_c}}\right)$ and $e_c+g_c<1$. Thus, applying the Chung's Lemma~\ref{lem:chung} to~\eqref{eq:c=1_1}, we obtain \begin{align}\label{eq:rate1}
    \E\|\bmu(t)-\ab_{t-1}^*\|^2 = O\left(\frac{1}{t^{h_c-g_c}}\right).
\end{align}
\textbf{Step 4. Convergence rate of $\E\|\bmu(t)-\ab^*\|^2$.}
Finally,~\eqref{eq:rate1} together with Lemma~\ref{lem:epsApprox} implies the result of theorem, namely, 
\begin{align*}
    \E\|\bmu(t)-\ab^*\|^2 =O\left(\frac{1}{t^{\min\{e_c,h_c-g_c\}}}\right), \, \mbox{for $c=1, 2$}. \qed
\end{align*}

\end{pf}

\subsection{Optimized convergence rate}
Defining the optimal settings for the parameters $\gamma_t$, $\varepsilon_t$, $\sigma_t$, and $\rho_t$ in Theorem~\ref{th:main}, we get the following result for the convergence rate.
\begin{corollary}
\label{cor:opt_rate}
Let Assumptions~\ref{assum:convex}-\ref{assum:infty} hold. Let in Algorithm~\ref{alg:algorithm1} the parameters be chosen  as follows:
\begin{align*}
    &\gamma_t = \frac{G_1}{t^{3/4+0.5\delta_1}}, \, \varepsilon_t=\frac{E_1}{t^{1/4-\delta_1}}, \, 
      \sigma_t=\frac{S_1}{t^{1/4}}, \, \rho_t =\frac{R_1}{t^{1/4(1-\delta_1)}},\cr
      &\hspace{6.2cm} \mbox{ if $c=1$;}\cr
    &\gamma_t = \frac{G_2}{t^{1/2}}, \, \varepsilon_t=\frac{E_2}{t^{1/2-\delta_2}}, \, 
      \sigma_t=\frac{S_2}{t}, \, \rho_t =\frac{R_2}{t^{1-0.5\delta_2}}, \cr
      &\hspace{6.2cm}  \mbox{ if $c=2$},
\end{align*}
where $\delta_1$ and $\delta_2$ are chosen to be arbitrary small. 
Moreover, let the joint action $\ab^*$ be a v-GNE of the game $\Gamma$. 
Then for the iterates Algorithm~\ref{alg:algorithm1} the following holds: 
\begin{align*}
    \E\|\bmu(t)-\ab^*\|^2 =\begin{cases}
        &O\left(\frac{1}{t^{1/4-\delta_1}}\right), \, \mbox{ if $c=1$},\\
        &O\left(\frac{1}{t^{1/2-\delta_2}}\right), \, \mbox{ if $c=2$}.
    \end{cases} 
\end{align*}
\end{corollary}
\begin{remark}
    Observe that the above result is the first convergence rate estimations in zeroth-order learning of GNE. Unfortunately (and naturally) it is worse than the best rates $O(1/\sqrt t)$ and  $O(1/t)$ attainable under the zeroth-order one- and two-point feedback in strongly monotone  games with uncoupled constraints~\cite{TatKamECC24}. 
   In comparison to works establishing rates under first-order information setting, ~\cite{jordan2023first} proves the rate $O(1/t^2)$ for the case of strongly monotone games with affine constraints,  but for convergence to a specifically defined surrogate of a GNE solution; while the best rate in strongly monotone games with uncoupled constraints is $O\left(\exp\{-\frac{t}{\gamma^2}\}\right)$ with $\gamma = L/\nu$, see~\cite{NesterovScrimali}. 
 
\end{remark}

One of the challenges in tightening the convergence rate in our setting isthe loose upper bound provided in Lemma~\ref{lem:epsApprox}. The next subsection discusses this issue and provides a setting where this bound can be tightened and, thus, the convergence rate estimation for Algorithm~\ref{alg:algorithm1} can be improved.
   
\subsection{The case of interior variational Nash equilibria}
\label{sub:th_intNE}
The analysis of the proposed algorithm in Theorem~\ref{th:main} is based on estimating the distance between the iterates $\bmu(t)$ and the primal part $\ab^*_t$ of the solution  to the regularized variational inequality~\eqref{eq:VI2}, see inequalities~\eqref{eq:c=1_1}, as well as the distance between $\ab^*_t$ and the solution $\ab^*$ (v-GNE) to the original game  (Lemma~\ref{lem:epsApprox}). The latter squared distance is proven not to exceed $\varepsilon_t\|\blab^*\|\|\blab^*-\blab^*_t\|$. Hence, to tighten this rate, one needs to get tight bounds on $\|\blab^*-\blab^*_t\|$.

Though $\|\blab^*-\blab^*_t\|$ converges to zero as $\varepsilon_t\to 0$ (see proof of Lemma~\ref{lem:epsApprox} in Appendix~\ref{app1}), it is generally not possible to provide\footnote{See discussion in Chapter 5.8 of the book~\cite{bakushinsky2012ill} on convergence rate properties of the regularized algorithms.} the order of the distance $\|\blab^*-\blab^*_t\|$. Thus, the provided order $\| \ab_t^*-\ab^*\|^2 = O(\varepsilon_t)$ is not tight.
This fact prevents improvement in the final rate estimation as the parameter $\varepsilon_t$ needs to be chosen to balance $\gamma_t$ ($\varepsilon_t + \gamma_t > 1/t$ as stated by Theorem~\ref{th:main}). 

By directly obtaining a tighter convergence rate for the sequence  $\|\ab^*_t-\ab^*\|^2$ one can prove a tighter rate estimation for the algorithm's iterates. The next statements demonstrate this fact in  the case for which the uncoupled game $\Gamma^{\blab}$ possesses an interior minimum norm Nash equilibrium. 

\begin{lem}\label{lem:lem_intNE}
Let the mapping $\Wb(\ab,\blab)$ be as in Definition~\ref{eq:gamemapping_aug} with $\Mb(\ab)$ satisfying Assumptions~\ref{assum:strmon} and~\ref{assum:infty}. 
    Let $[\ab^*,\blab^*]$ be a minimal norm solution to~\eqref{eq:VI1} and $[\ab^*_t,\blab^*_t]$ be a solution to~\eqref{eq:VI2}. Moreover, let $\ab^*$ belong to the interior of the set $\Ab$. Then, as $\varepsilon_t\to 0$, $[\ab^*_t,\blab^*_t]$ converges to $[\ab^*,\blab^*]$ and, moreover,
    \[\|\ab^*-\ab^*_t\|\le \varepsilon_t\frac{\|\blab^*\|L}{\|K\|\nu},\]
    where $L$ is the Lipschitz constant of the pseudo-gradient $\Mb$ (see Remark~\ref{rem:Lip}).
\end{lem}
\begin{pf}
According to Definition~\ref{eq:gamemapping_aug} and~\eqref{eq:primalM} and due to the fact that $[\ab^*,\blab^*]$ solves~\eqref{eq:VI1} and $\ab^*$ belongs to the interior of $\Ab$, we conclude that
$\Wb^{pr}(\ab^*, \blab^*) = \Mb(\ab^*) + K^T\blab^* = \boldsymbol{0}$.
On the other hand, according to the result in Lemma~\ref{lem:epsApprox}, for the distance between $\ab^*$ and the primal part $\ab^*_t$ of the solution $\zb^*_t = [\ab^*_t,\lambda_t^*]$ to~\eqref{eq:VI2} the following relation holds: 
$\lim_{\varepsilon_t\to 0}\|\ab^*_t - \ab^*\| = 0$.
Thus, there exists a sufficient large $T$ such that the solution $\ab^*_t$ to~\eqref{eq:VI2}, $t\ge T$, belongs to the interior of $\Ab$ as well and, thus,  
$\Wb^{pr}(\ab_t^*,\blab_t^*) = \Mb(\ab_t^*) + K^T \blab_t^* = \boldsymbol{0}$.
Hence, 
$\Mb(\ab^*) + K^T\blab^* = \Mb(\ab_t^*) + K^T\blab_t^*$, which implies that 
\begin{align*}
    \|\blab^*-\blab_t^*\|
    \|K\| = \|\Mb(\ab^*) - \Mb(\ab_t^*)\| \le L\|\ab^*-\ab_t^*\|,
\end{align*}
where the last inequality is due to Assumption~\ref{assum:infty} (see Remark~\ref{rem:Lip}). Next, using Lemma~\ref{lem:epsApprox}, we conclude that 
\[\|\ab^*-\ab^*_t\|^2 \le\varepsilon_t\frac{\|\blab^*\|\|\blab^*-\blab^*_t\|}{\nu}\le\varepsilon_t\frac{\|\blab^*\|L\|\ab^*-\ab_t^*\|}{\|K\|\nu}.\]
Thus, 
$\|\ab^*-\ab^*_t\|\le \varepsilon_t\frac{\|\blab^*\|L}{\|K\|\nu}$.
$\qed$\end{pf}
Having the lemma above in place, we obtain the following result on the convergence rate for the special case of the interior Nash equilibrium. 

\begin{theorem}\label{th:th_intNE}
    Let the conditions of Theorem~\ref{th:main} hold. Moreover, let the unique v-GNE in the game $\Gamma$ belong to the interior of the set $\Ab$. Then for the iterates in Algorithm~\ref{alg:algorithm1} the following estimation holds: 
\begin{align*}
    \E\|\bmu(t)-\ab^*\|^2 =O\left(\frac{1}{t^{\min\{2e_c,h_c-g_c\}}}\right), \, \mbox{for $c=1, 2$}.
\end{align*}
\end{theorem}
\begin{pf}
Taking Lemma~\ref{lem:lem_intNE} into account, we adjust \textbf{Part~4} in the proof of Theorem~\ref{th:main} as follows:
\begin{align*}
    \E\|\bmu(t)-\ab^*\|^2 =O\left(\frac{1}{t^{\min\{2e_c,h_c-g_c\}}}\right), \, \mbox{$c=1, 2$}. \qed
\end{align*}
\end{pf}
\begin{corollary}\label{cor:int}
Optimizing the parameters $e_c$, $s_c$, $r_c$, and $g_c$, we obtain the following convergence rate under the conditions of Theorem~\ref{th:th_intNE} 
\begin{align*}
    \E\|\bmu(t)-\ab^*\|^2 =\begin{cases}
        &O\left(\frac{1}{t^{4/15-\delta_1}}\right), \, \mbox{ if $c=1$},\\
        &O\left(\frac{1}{t^{4/7}}\right), \, \mbox{ if $c=2$}.
    \end{cases} 
\end{align*} 
for any small $\delta_1,\delta_2>0$ under the choice $e_1 = 2/15$, $s_1 = 4/15$, $r_1 = 4/15-\delta_1$, $g_1 = 4/5$ and $e_2=2/7,$ $r_2\ge 2/7$, $s_2>r_2$, $g_2 = 4/7$ respectively. 
\end{corollary}


Note that the interior condition in Theorem~\ref{th:th_intNE} holds, whenever $\Ab = \R^D$ and solution to $\Gamma^{\blab}$ exists. However, in   this case compactness of $\Ab$ is not fulfilled. Thus, to be able to repeat the proof of Theorem~\ref{th:main} one will need a preliminary result guaranteeing boundedness of the iterates $\bmu(t)$ in Algorithm~\ref{alg:algorithm1}.  The work~\cite{ECC25} provides the corresponding analysis in the case of two-point estimates ($c=2$) and obtains the rate of the order $O\left(\frac{1}{t^{4/7}}\right)$ which is in consistence with Corollary~\ref{cor:int}.
Extending the analysis to the case of one-point gradient estimates (the case $c=1$) requires further assumptions to ensure boundedness of iterates and is a potential direction for further investigations. 

\section{Numerical example}
\label{sec:example}

\begin{figure*}[!h]
	\centering
	\begin{minipage}[t]{0.45\textwidth}
		\centering
		\includegraphics[width=\textwidth]{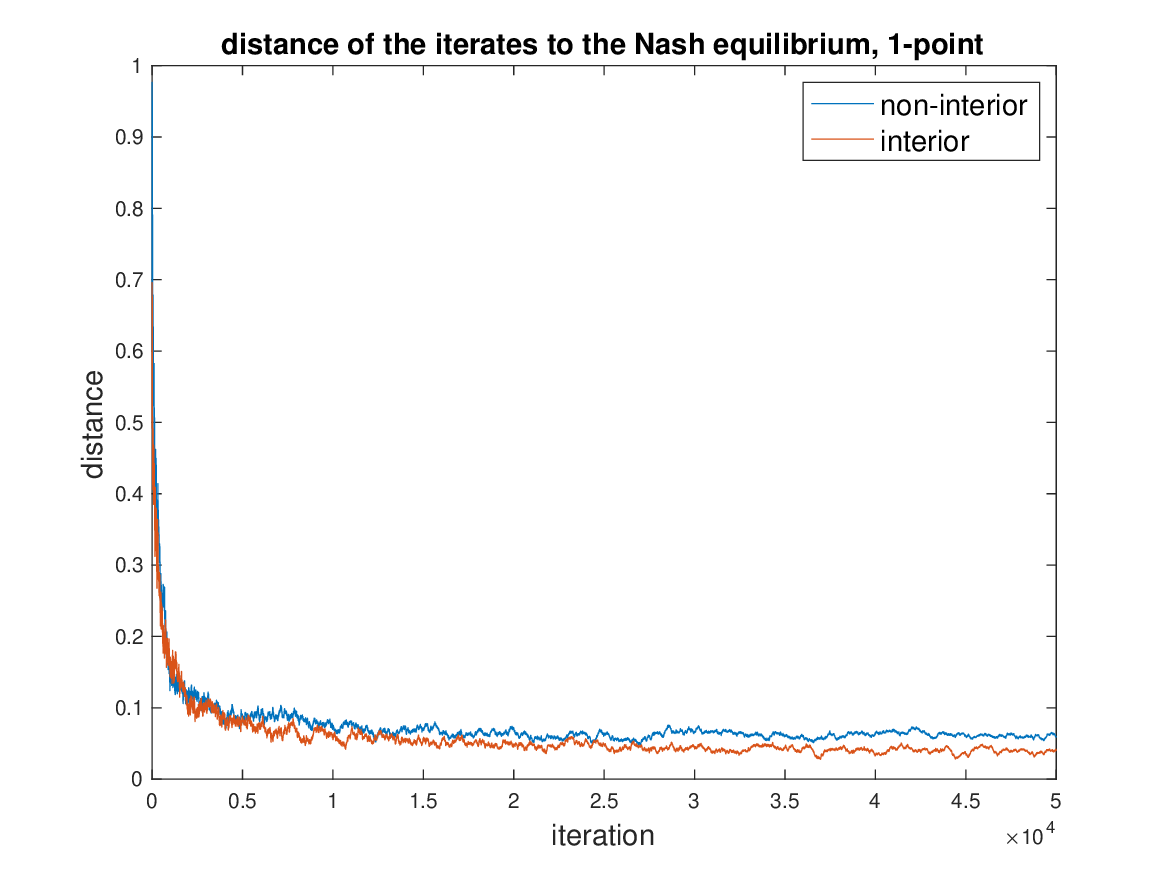}
		\caption*{(a) One-point feedback}
		\label{fig:d1point}
	\end{minipage}
	\hfill
	\begin{minipage}[t]{0.45\textwidth}
		\centering
		\includegraphics[width=\textwidth]{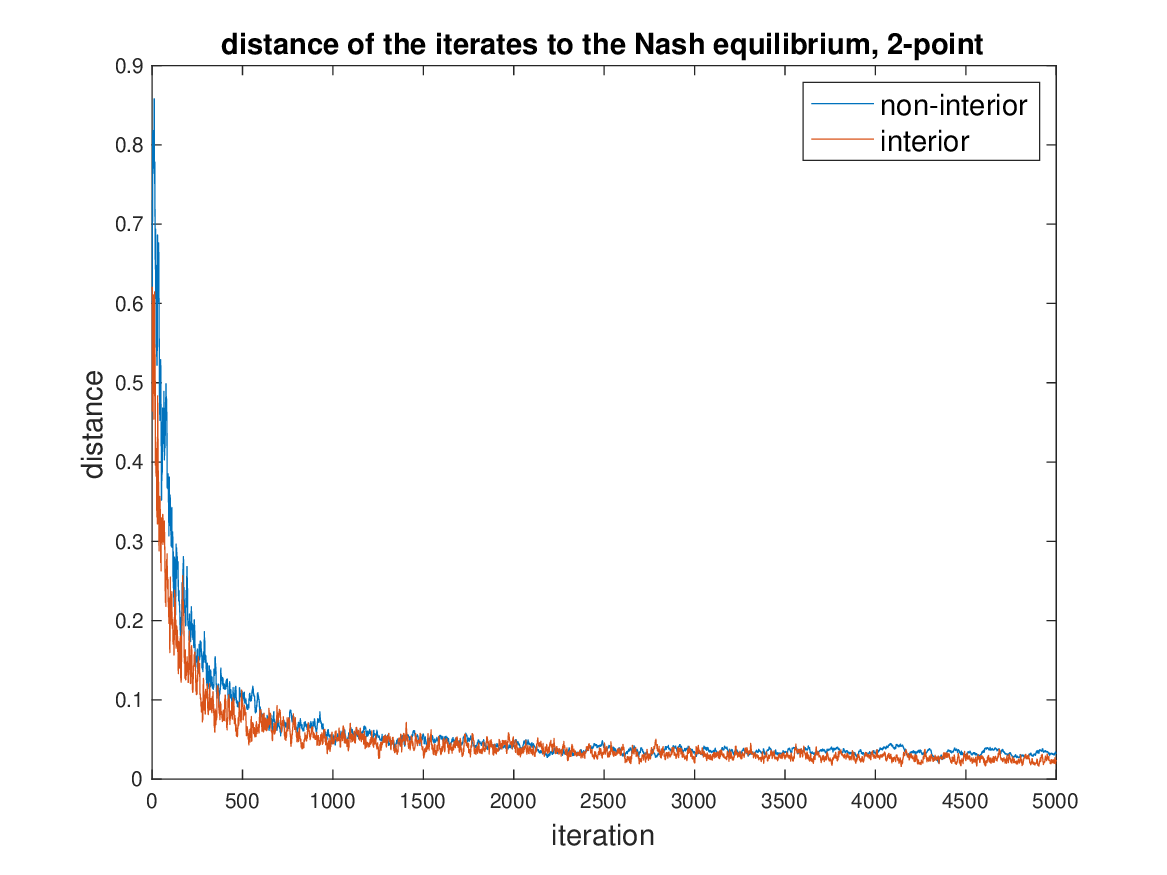}
		\caption*{(b) Two-point feedback}
		\label{fig:d2point}
	\end{minipage}
	\caption{Average distance of the iterates from the Nash equilibrium in case of one-point (left) and two-point (right) feedback. The average is computed over 20 runs of Algorithm~\ref{alg:algorithm1}, each starting from the same initial condition for $\bmu^i$, $i=1, \dots, N$ and $\bla$.}
	\label{fig:distance}
\end{figure*}

Consider a discrete-time dynamical system with two states, which is steered by two agents. Denote the states at time $t$ by  $s_{1,t}, s_{2,t} \in \R$ and the input of player $i$ to state $j$ at time $t$ by $a^i_{j,t}$ and the state of the system at time $t$ by $s_{j,t}$. Given an initial condition of $s_{1,0}, s_{2,0} \in \R$, the state of the system evolves over time as
\begin{align}
\label{eq:dynamics}
s_{1,t+1} &= f_1(s_{1,t},a^1_{1,t}, a^2_{1,t}), \\
\nonumber
s_{2,t+1} &= f_2(s_{2,t}, a^1_{2,t}, a^2_{2,t}).   
\end{align}
Each agent has a preference on the states $(s_1, s_2)$ being close to a nominal value  $(\bar{s}^i_1, \bar{s}^i_2)$ at the end of a given horizon $T$. Thus, the agents' objectives are coupled which leads to a game between the agents/players. As for concrete applications of the setup, consider two players steering a car with different objectives or, see an application in competitive supply-chain in \cite{hall2024receding}.

Each agent also has  a cost corresponding to their inputs. Let $\ab = (\ab^1, \ab^2) \in \R^{4T}$, with \[\ab^i = (a^i_{1,0}, a^i_{2,0}, \dots, a^i_{1,T-1},a^i_{2,T-1}) \in \R^{2T},\] for $i \in [1,2]$. The following  cost function  captures  player $i$'s objective:
\begin{align*}
&J^i(\ab) = \frac{1}{2}  \big[ \begin{pmatrix} s_{1,T}-\bar{s}^i_{1},\\ s_{2,T}-\bar{s}^i_2\end{pmatrix}^TQ^i\begin{pmatrix} s_{1,T}-\bar{s}^i_{1},\\ s_{2,T}-\bar{s}^i_2\end{pmatrix}  \\
+& \sum_{t=0}^{T-1}   \begin{pmatrix} a^i_{1,t} \\ a^i_{2,t}\end{pmatrix}^T R^i \begin{pmatrix} a^i_{1,t} \\ a^i_{2,t}\end{pmatrix} + h^i \begin{pmatrix} a^i_{1,t} \\ a^i_{2,t}\end{pmatrix}\big],
\end{align*}
where $Q^i, R^i \in \R^{2 \times 2}$, for $i=1,2$ are positive definite matrices and $h^i \in \R^2$ is a weight vector. 

There are bounds on the magnitude of the inputs of each player and  on the magnitude of input to each state. Thus, for $i=1,2$ we have $0 \leq a^i_{j,t} \leq \bar{a}^i$, for each of the inputs $j=1,2$ and at each time $t =0, 1, \dots, T-1$, where $\bar{a}^i$ denote the maximum input  each player can apply to each state. The bound on total input leads to a coupling constraint given by $| a^1_{j,t} + a^2_{j,t} | \leq c_j$, for $j=1,2$, where $c_j$ is the maximum allowable input to each state.

Due to the coupling in the cost and constraints, the problem falls into the class of GNE problems with jointly convex constraints. The goal  for each player is to compute a v-GNE of the game, without knowledge of the dynamics $f_1, f_2$ and others' choice of actions. However, with a payoff-based information, she can observe her cost function for any jointly chosen action $\ab \in \R^{2T}$. 

The monotonicity of the game depends on the specific model of the dynamics $f_1, f_2$ in \eqref{eq:dynamics} and the interplay of the matrices $Q^i, R^i$ for $i=1, 2$. For concreteness, let us consider the following simple dynamics in a time horizon of one-step,  and drop the $t$ from the subindex for simplicity in notation: $s_1 = s_{1,0} + a^1_1 + a^2_1$, $s_2 = s_{2,0} + a^1_2 + a^2_2$.
Let the cost be 
\begin{align*}
&J^i(\ab) = \frac{1}{2}  \big[ q_1^i (s_{1}-\bar{s}^i_{1})^2 + q_2^i (s_{2}-\bar{s}^i_2)^2 \\
+& r^i (a^i_{1})^2 + r^i  (a^i_{2})^2 + r^i(a^i_1 + a^i_2 - \bar{a}^1)^2\big],
\end{align*}
It can be verified that \[\Mb(\ab) = [\frac{\partial J_1(\ab)}{\partial a_1^1}, \frac{\partial J_1(\ab)}{\partial a_2^1}, \frac{\partial J_2(\ab)}{\partial a_1^2}, \frac{\partial J_2(\ab)}{\partial a_2^2}]^T\] is strongly monotone. Furthermore, since $\nabla \Mb$ is non-symmetric, the game is not potential. 

We apply our algorithm to compute the Nash equilibrium. The parameters are set as follows:  $s_{1,0} = s_{2,0} = 1/3$, $(\bar{s}^1_1, \bar{s}^1_2,\bar{s}^2_1,\bar{s}^2_2) = (1, 2/3, 2/3, 2/3)$, $(q^1_1, q^1_2, q^2_1, q^2_2) = (5,3,1,5)$, and $r^1 = 2, r^2 =1$. The coupling constraint is set to $(\bar{c}_{1}, \bar{c}_{2}) = (2/3,5/3)$.  Two cases for the individual constraints are considered. In case 1, they are set to $\bar{a}^1 = \bar{a}^2 =1 $ and in case 2, $\bar{a}^1 = 0.3$, $\bar{a}^2 =1 $. In the first case,  we find that the solution of ${\Mb}(\ab) = 0$, given as $(0.5246, 0.0352, 0.1252, 0.4332)$, satisfies all individual and coupling constraints and thus, is a Nash equilibrium. Thus, we have an interior Nash equilibrium. In the second case, the solution found above does not satisfy the constraint on $\bar{a}^1$. However, by solving  the KKT system of equations centrally using the first order information, we find the Nash equilibrium to be $(0.3, 0, 0.1950, 0.4483)$.

Figure \ref{fig:distance} shows the distance of the iterates to the v-GNE, $\| \bmu(t) - \ab^*\|$ in the case of one-point (left) and two-point (right) feedback. Notice that in both cases, two-point feedback significantly speeds up the convergence as the number of iterates to approach an equilibrium has decreased by an order of magnitude (from $10^4$ to $10^3$). Moreover, from these figures, we can also see a small improvement in convergence speed for the case of the interior Nash equilibrium (red) versus the non-interior equilibrium (blue).

\section{Conclusion}

In this work, we addressed the problem of payoff-based learning in generalized Nash equilibrium (GNE) problems with linear coupling constraints. We proposed a novel algorithm that, under the assumptions of strongly monotone and Lipschitz-continuous pseudo-gradients, achieves convergence to a variational GNE (v-GNE). The algorithm handles both one-point and two-point feedback settings, and relies on Tikhonov regularization to overcome the challenges caused by the lack of strong monotonicity in the extended pseudo-gradient. We rigorously analyzed the convergence behavior of the algorithm and established explicit non-asymptotic convergence rates in both settings. To our knowledge, this is the first work to derive such rates in the payoff-based setting for games with affine constraints and without potential structure.
We demonstrated that the proposed algorithm achieves a convergence rate of $O(1/t^{1/4 - \delta})$ in the one-point feedback setting and $O(1/t^{1/2 - \delta})$ in the two-point feedback setting, with potential improvement in the presence of interior solutions.  Future work may explore extensions to games with player-wise coupling constraints. Furthermore, it would be important to establish lower bounds on convergence rate of payoff-based learning for the considered class of games.


\bibliographystyle{plain}        
\bibliography{autosam}           



\appendix
\section{Proof of Lemma~\ref{lem:exist_uncoupled}}\label{app:exist_uncoupled}
\begin{pf} According to the definition of variational GNE (see Definition~\ref{def:VGNE}), $\ab^*$ is variational GNE iif it solves~\eqref{eq:GVI}.  
Thus, the assertion 1) holds due to Proposition 1.3.4. in~\cite{FaccPang1}, as a constraint qualification condition is automatically fulfilled in the case of affine constraints under consideration (the so called relaxed Slater's condition implying no duality gap in the convex case, see~\cite{Lecture_Nemirovski}). 

Existence of a unique variational GNE in $\Gamma$ is guaranteed under Assumption~\ref{assum:strmon}, as there exists a unique solution to~\eqref{eq:GVI}, if the mapping $\Mb$ is strongly monotone. Thus, the assertion 3) is implied. Moreover, given Assumption~\ref{assum:convex}  and taking into account 1), the assertion 2) follows. 
Finally, the assertion 4) holds, as demonstrated in Lemma~5.1 in \cite{ZhuFrazzoli}.
$\qed$\end{pf}

\section{Proof of Lemmas~\ref{lem:epsApprox} and~\ref{lem:t_vs_t-1}}\label{app1}
(Proof of Lemma~\ref{lem:epsApprox}.)
    According to the definition of $\Wb_t$ and  given $\zb^* = [\ab^*, \blab^*]$, $\zb^*_t =  [\ab^*_t, \blab^*_t]$, we have
    \begin{align*}
        &\langle \Wb(\zb^*),\zb^*_t-\zb^*\rangle\ge 0,\cr
        &\langle \Wb_t(\zb_t^*),\zb^*-\zb_t^*\rangle=\langle \Wb(\zb_t^*)+\varepsilon_tI_{0,n}\zb_t^*,\zb^*-\zb_t^*\rangle\ge 0.
    \end{align*}
    By summing up the inequalities above and using~\eqref{eq:aug-psedo} (holds due to Assumption~\ref{assum:strmon}), we obtain
    \begin{align}\label{eq:sum}
        0\le -\nu\|\ab^*-\ab^*_t\|^2 +\varepsilon_t\langle\blab^*_t,\blab^*-\blab^*_t\rangle,
    \end{align}
    which implies
    \begin{align*}
        \nu\|\ab^*-\ab^*_t\|^2 +\varepsilon_t\|\blab^*-\blab^*_t\|^2\le\varepsilon_t\langle\blab^*\blab^*-\blab^*_t\rangle.
    \end{align*}
    Thus, on the one hand, $\varepsilon_t\|\blab^*-\blab^*_t\|^2\le\varepsilon_t\|\blab^*\|\|\blab^*-\blab^*_t\|$,
       implying 
   $\|\blab^*-\blab^*_t\|\le\|\blab^*\|$ and, thus, uniform boundedness of $\|\blab^*-\blab^*_t\|$. Moreover, as  $\blab^* = \Proj_{\R^n_{+}}[\blab^* + \alpha \gb(\ab^*)]$, $\blab_t^* = \Proj_{\R^n_{+}}[\blab_t^* + \alpha (\gb(\ab_t^*) + \varepsilon_t\blab_t^*)]$ for any $\alpha\ge 0$, and the projection is a continuous mapping, we conclude that  $\lim_{\varepsilon_t\to 0}\|\blab^*-\blab^*_t\| = 0$, where $[\ab^*,\blab^*]$ is  the minimal norm solution to~\eqref{eq:VI1}. On the other hand,~\eqref{eq:sum} together with the previous relation implies
   \begin{align*}
        \nu\|\ab^*-\ab^*_t\|^2 \le\varepsilon_t\|\blab^*\|\|\blab^*-\blab^*_t\|\le\varepsilon_t\|\blab^*\|^2.   
    \end{align*}

(Proof of Lemma~\ref{lem:t_vs_t-1}.) According to 
the definition of $\Wb_t$ and  given $\zb^*_t = [\ab^*, \blab^*]$, we have
    \begin{align*}
        &\langle \Wb(\zb_{t-1}^*)+\varepsilon_{t-1}I_{0,n}\zb_{t-1}^*,\zb^*_t-\zb_{t-1}^*\rangle\ge 0,\cr
        &\langle \Wb(\zb_t^*)+\varepsilon_tI_{0,n}\zb_t^*,\zb_{t-1}^*-\zb_t^*\rangle\ge 0.
    \end{align*}
    By summing up the inequalities above and using~\eqref{eq:aug-psedo}, we obtain
    \begin{align}\label{eq:sum2}
        &0\le-\nu\|\ab^*-\ab^*_t\|^2 -\varepsilon_t\|\blab^*-\blab^*_t\|^2 \cr
        &\qquad+ (\varepsilon_{t-1}-\varepsilon_t)\langle \blab_{t-1}^*,\blab^*_t-\blab_{t-1}^*\rangle,
    \end{align}
    which implies, on the one hand,
    \begin{align*}
        &\varepsilon_t\|\blab^*-\blab^*_t\|^2 \le |\varepsilon_{t-1}-\varepsilon_t| \|\blab_{t-1}^*\|\|\blab^*_t-\blab_{t-1}^*\|.
    \end{align*}
    Thus, due to the boundedness of $\|\blab_{t-1}^*\|$ (see the inequality $\|\blab^*-\blab^*_t\|\le\|\blab^*\|$ in the proof of Lemma~\ref{lem:epsApprox}), we conclude that
    \[\|\blab^*_t-\blab^*_{t-1}\|^2 = O\left(\frac{(\varepsilon_t-\varepsilon_{t-1})^2}{\varepsilon^2_{t}}\right).\]
    On the other hand, the relation above together with~\eqref{eq:sum2}, implies
    \[\|\ab^*_t-\ab^*_{t-1}\|^2 = O\left(\frac{(\varepsilon_t-\varepsilon_{t-1})^2}{\varepsilon_{t}}\right). \]

\section{Theorem 2.5.2 in \cite{NH}}\label{app_bound}
\begin{theorem}\label{th:bound}
Let $V(t,X)$ be a rel-valued function defined for some  stochastic process $\{X(t)\}_t$ taking values in $\R^d$. Let us assume existence of the following generating operator 
\begin{align*}
LV(t,x)=E[V(t+1, X(t+1))\mid X(t)=x]-V(t,x),
\end{align*}
for which the following decay holds: 
\begin{align*}
LV(t,x)\le -\alpha(t+1)\psi(x) + \phi(t)(1+V(t,x)),
\end{align*}
where $\psi\ge 0$ on $ \R^{d}$, $\phi(t)>0$, $\forall t$, $\sum_{t=0}^{\infty}\phi(t)<\infty$, $\alpha(t)>0$, $\sum_{t=0}^{\infty} \alpha(t)= \infty$. Then the process $X(t)$ is almost sure boundeded over $t$.
\end{theorem}

\section{Chung's Lemma}\label{app2}
Here we formulate the result from Lemma 4 in \cite{Chung}.
\begin{lem}\label{lem:chung}
Suppose that $\{b_n\}$, $n\ge1$, is a sequence of real numbers such that for $n\ge n_0$,
\[b_{n+1}\le\left(1-\frac{c_n}{n^s}\right)b_n+\frac{c'}{n^t},\]
where $0<s<1$, $c_n\ge c>0$, $c'>0$. Then
\[\limsup_{n\to\infty}n^{t-s}b_n\le \frac{c'}{c}.\]
\end{lem}
\end{document}